\documentclass[11pt,reqno]{amsart}
\usepackage{amsfonts,mathrsfs,mathabx,bezier}
\usepackage{color}
\usepackage{makecell,multirow,diagbox}
\usepackage{amsmath,amsxtra,latexsym,amsthm,amssymb,amscd,pb-diagram}
\usepackage[colorlinks=true,linkcolor=magenta,citecolor=blue]{hyperref}
\usepackage[margin=1in]{geometry}

\usepackage{graphicx}
\usepackage{subfigure}
\usepackage{enumitem}

\newtheorem{dn}{Definition}[section]
\newtheorem{dl}{Theorem}[section]
\newtheorem{md}{Proposition}[section]
\newtheorem{bd}{Lemma}[section]
\newtheorem{hq}{Corollary}[section]
\newtheorem{nx}{Remark}[section]
\newtheorem{vd}{Example}[section]
\newcommand{\R}{\mathbb{R}}

\newcommand{\Z}{\mathbb{Z}}
\newcommand{\N}{\mathbb{N}}
\newcommand{\e}{\varepsilon}
\newcommand{\ity}{\infty}

\newcommand{\bbd}{\begin{bd}}
\newcommand{\ebd}{\end{bd}}
\newcommand{\bdn}{\begin{dn}}
\newcommand{\edn}{\end{dn}}
\newcommand{\bhq}{\begin{hq}}
\newcommand{\ehq}{\end{hq}}
\newcommand{\bdl}{\begin{dl}}
\newcommand{\edl}{\end{dl}}
\newcommand{\bnx}{\begin{nx}}
\newcommand{\enx}{\end{nx}}
\newcommand{\bmd}{\begin{md}}
\newcommand{\emd}{\end{md}}
\newcommand{\bvd}{\begin{vd}}
\newcommand{\evd}{\end{vd}}

\title[Weakly coupled system of damped wave equations with moduli of continuity]{The interplay of critical regularity of nonlinearities in a weakly coupled system of semi-linear damped wave equations}
\author{Tuan Anh Dao}
\address{Tuan Anh Dao \hfill\break
$\quad$ School of Applied Mathematics and Informatics, Hanoi University of Science and Technology, No.1 Dai Co Viet road, Hanoi, Vietnam \hfill\break
Faculty for Mathematics and Computer Science, TU Bergakademie Freiberg, Pr\"{u}ferstr. 9, 09596, Freiberg, Germany}
\email{anh.daotuan@hust.edu.vn}
\author{Michael Reissig}
\address{Michael Reissig \hfill\break
$\quad$ Faculty for Mathematics and Computer Science, TU Bergakademie Freiberg, Pr\"{u}ferstr. 9, 09596, Freiberg, Germany}
\email{reissig@math.tu-freiberg.de}

\begin{document}
\subjclass[2010]{26A15, 35A01, 35L52}
\keywords{Damped wave equations; Weakly coupled system; Modulus of continuity; Global existence; Blow-up}
% \thanks{$^*$ Corresponding author}

\begin{abstract}
We would like to study a weakly coupled system of semi-linear classical damped wave equations with moduli of continuity in nonlinearities whose powers belong to the critical curve in the $p-q$ plane. The main goal of this paper is to find out the sharp conditions of these moduli of continuity which classify between global (in time) existence of small data solutions and finite time blow-up of solutions.
\end{abstract}

\maketitle

%==========================================================Introduction
\section{Introduction}
At first, concerning the following classical semi-linear damped wave equation:
\begin{equation} \label{pt1.1}
\begin{cases}
u_{tt}- \Delta u+ u_t= |u|^p,& \quad x \in \R^n,\, t\ge 0, \\
u(0,x)= u_0(x), \quad u_t(0,x)= u_1(x),& \quad x \in \R^n,
\end{cases}
\end{equation}
with $p>1$, the authors in \cite{TodorovaYordanov} proved the global (in time) existence of energy solutions for
$$p>p_{\text{Fuj}}(n)= 1+ \frac{2}{n}, $$
the so-called Fujita exponent, and for $p\le n/(n-2)$ if $n\ge 3$. Besides, they also indicated a blow-up result in the inverse case $1<p<p_{\text{Fuj}}(n)$ which was improved for $1<p\le p_{\text{Fuj}}(n)$ in the paper \cite{Zhang} by using the well-known test function method so far. For this reason, we can say that the Fujita exponent distinguishes the admissible range of powers $p$ in (\ref{pt1.1}) into those possessing global existence of small data solutions (stability of zero solutions) and those producing a blow-up behavior. However, to determine the critical exponent $p_{\text{Fuj}}(n)$, it seems too rough to restrict (\ref{pt1.1}) in the scale of power nonlinearities $\{|u|^p\}_{p>1}$. Quite recently, the second author and his collaborators have discussed this issue for the following Cauchy problem in the paper \cite{EbertGirardiReissig}:
\begin{equation} \label{pt1.2}
\begin{cases}
u_{tt}- \Delta u+ u_t= |u|^{p_{\text{Fuj}}(n)}\mu(|u|),& \quad x \in \R^n,\, t\ge 0, \\
u(0,x)= u_0(x), \quad u_t(0,x)= u_1(x),& \quad x \in \R^n,
\end{cases}
\end{equation}
where $\mu= \mu(|u|)$ stands for a modulus of continuity, a well-known notation to describe the regularity of a function with respect to desired variables. This means that it provides an additional regularity of the nonlinear term in comparison with the power nonlinearity $|u|^{p_{\text{Fuj}}(n)}$. More precisely, the authors in the cited paper have found out the sharp conditions for the critical regularity of the nonlinear term of (\ref{pt1.1}), namely,
$$ \int_0^c \frac{\mu(s)}{s}ds< \ity \quad \text{ and }\quad \int_0^c \frac{\mu(s)}{s}ds= \ity, $$
where $c$ is a sufficiently small positive constant, which separates the global (in time) existence of small data Sobolev solutions and the blow-up behavior of Sobolev solutions, respectively.

During the last decades, the study of Cauchy problems for weakly coupled systems of equations in place of exploiting single equations only has been achieving a great attention from many mathematicians because of their wide applications in various disciplines. One of the most typical problems is the following weakly coupled system of semi-linear classical damped wave equations (see, for example, \cite{Narazaki,NishiharaWakasugi,SunWang,Takeda}):
\begin{equation} \label{pt1.3}
\begin{cases}
u_{tt}- \Delta u+ u_t= |v|^p, &\quad x\in \R^n,\, t \ge 0, \\
v_{tt}- \Delta v+ v_t= |u|^q, &\quad x\in \R^n,\, t \ge 0, \\
u(0,x)= u_0(x),\quad u_t(0,x)=u_1(x), &\quad x\in \R^n, \\
v(0,x)= v_0(x),\quad v_t(0,x)=v_1(x), &\quad x\in \R^n,
\end{cases}
\end{equation}
with $p,q >1$. Particularly, the authors in the papers \cite{Narazaki,SunWang} investigated the global (in time) existence of small data solutions to (\ref{pt1.3}) in low space dimensions $n=1,2,3$, which was extended for any space dimensions $n\ge 1$ in the paper \cite{NishiharaWakasugi} afterwards by using weighted energy estimates, if the following condition for a pair $(p,q)$ holds:
$$ \frac{1+ \max\{p,q\}}{pq-1}< \frac{n}{2}. $$
Additionally, when the above condition is no longer true, non-existence results of global (in time) solutions to (\ref{pt1.3}) were proved in \cite{Narazaki,NishiharaWakasugi,SunWang}. For this reason, one can claim that the critical curve of the power exponents for (\ref{pt1.3}) in the $p-q$ plane is
\begin{equation}
\frac{1+ \max\{p,q\}}{pq-1}= \frac{n}{2}. \label{pt1.5}
\end{equation}

Our interest of this paper is strongly inspired by the recent paper of the second author \cite{EbertGirardiReissig} in the connection between (\ref{pt1.2}) and (\ref{pt1.3}). A natural question arises that whether it is sharp or not to obtain the critical curve (\ref{pt1.5}) in the scale of pairs of power nonlinearities $\{|v|^p,|u|^q\}_{p,q>1}$. Hence, the key motivation for this article is to give a positive answer to this question. Namely, let us consider the Cauchy problem for weakly coupled system of semi-linear damped wave equations with moduli of continuity term in nonlinearities as follows:
\begin{equation} \label{pt1.4}
\begin{cases}
u_{tt}- \Delta u+ u_t= |v|^{p_c}\,\mu_1(|v|), &\quad x\in \R^n,\, t \ge 0, \\
v_{tt}- \Delta v+ v_t= |u|^{q_c}\,\mu_2(|u|), &\quad x\in \R^n,\, t \ge 0, \\
u(0,x)= u_0(x),\quad u_t(0,x)=u_1(x), &\quad x\in \R^n, \\
v(0,x)= v_0(x),\quad v_t(0,x)=v_1(x), &\quad x\in \R^n,
\end{cases}
\end{equation}
where the functions $\mu_1= \mu_1(|v|)$ and $\mu_2= \mu_2(|u|)$ are some suitable moduli of continuity. We assume that the powers $p_c,q_c >1$ of nonlinearities belong to the critical curve (\ref{pt1.5}) in the $p-q$ plane. Our main purpose of this paper is that we would like to understand the effect each other of additional regularities of power nonlinearities $|v|^{p_c}$ and $|u|^{q_c}$, which are given by these moduli of continuity, on not only global (in time) existence of small data solutions but also finite time blow-up of solutions. Especially, we are interested in looking for a threshold by exploring the following optimal conditions for $\mu_1$ and $\mu_2$:
$$ \int_0^{c}\frac{1}{s}\big(\mu_1(s)\big)^{\frac{q_c}{q_c+1}}\big(\mu_2(s)\big)^{\frac{1}{q_c+1}}\,ds<\ity \quad \text{ or }\quad \int_0^{c}\frac{1}{s}\big(\mu_1(s)\big)^{\frac{q_c}{q_c+1}}\big(\mu_2(s)\big)^{\frac{1}{q_c+1}}\,ds=\ity, $$
which leads to either global existence results or non-existence results of global solutions individually (see later, Theorems \ref{GlobalExistence1.4} and \ref{dloptimal1}), with developing some ideas coming from the paper \cite{EbertGirardiReissig} to the weakly coupled system of type (\ref{pt1.4}). Through this work, one should recognize that our results are not a simple generalization of those in \cite{EbertGirardiReissig}. Concretely, there are two points worthy to be mentioned. The first point as we can see is that allowing loss of decay appropriately, which has never appeared in \cite{EbertGirardiReissig}, comes into play to find out these conditions for $\mu_1$ and $\mu_2$ in guaranteeing the existence of global solutions. In other words, we can feel more explicit how the required assumptions of additional regularities of nonlinearities follow essentially from using some suitable loss of decay. This gives a new interplay in comparison with the previous research manuscripts in terms of the study of weakly coupled systems (see, for instance, \cite{Narazaki,NishiharaWakasugi,SunWang}). Moreover, the other point worth noticing is that the technical choice of a test function with the parameter depending on $p_c,q_c$ brings some remarkable benefits in the proof of blow-up part. \medskip

\textbf{Notations}
\begin{itemize}[leftmargin=*]
\item We denote $[s]:= \max \big\{k \in \Z \,\, : \,\, k\le s \big\}$ as the integer part of $s \in \R$.
\item For later convenience, hereafter $C$ denotes a suitable positive constant and may have different value from line to line.
\item For two given nonnegative functions $f$ and $g$, we write $f\lesssim g$ when $f\le Cg$. We write $f \approx g$ when $g\lesssim f\lesssim g$.
\item As usual, $H^m$ and $\dot{H}^m$, with $m \in \N$, denote Sobolev spaces based on the $L^2$ spaces.
\end{itemize}

\textbf{Main results} \medskip

Without loss of generality, if we assume $p_c<q_c$, then the critical curve in the $p-q$ plane for (\ref{pt1.4}) becomes
\begin{equation}
\frac{1+ q_c}{p_c q_c-1}= \frac{n}{2}. \label{CriticalCurve}
\end{equation}
Our main results concerned with the case $p_c<q_c$ read as follows.

\bdl[\textbf{Global existence}] \label{GlobalExistence1.4}
Let $n=1,2$. Assume that the following assumptions of moduli of continuity hold:
\begin{equation}
s\mu_j'(s) \lesssim \mu_j(s) \quad \text{ with }j=1,2. \label{*}
\end{equation}
Moreover, we suppose that one of the following conditions is satisfied:
\begin{align}
\text{\fontshape{n}\selectfont i)}&\,\, \int_0^{c}\frac{\mu_1(s)}{s}ds<\ity \quad \text{ and }\quad \int_0^{c}\frac{\mu_2(s)}{s}ds<\ity. \label{existenceCond.1.1.4} \\
\nonumber \\
\text{\fontshape{n}\selectfont ii)}&\,\, \text{ If }\,\, \int_0^{c}\frac{\mu_1(s)}{s}ds= \ity\,\, \text{ or }\,\, \int_0^{c}\frac{\mu_2(s)}{s}ds= \ity\,\,\text{ then }\,\, \int_0^{c}\frac{1}{s}\big(\mu_1(s)\big)^{\frac{q_c}{q_c+1}}\big(\mu_2(s)\big)^{\frac{1}{q_c+1}}\,ds<\ity. \label{existenceCond.2.1.4}
\end{align}
Here $c >0$ is a suitable small constant. Then, there exists a constant $\epsilon>0$ such that for any small data
$$ \big((u_0,u_1),\, (v_0,v_1) \big) \in \mathcal{A}:= \Big(\big(L^1 \cap H^{1+[n/2]} \big) \times \big(L^1 \cap H^{[n/2]}\big)\Big)^2 $$
satisfying the assumption
$$ \|(u_0,u_1),(v_0,v_1)\|_{\mathcal{A}}:= \|u_0\|_{H^{1+[n/2]}}+ \|u_0\|_{L^1}+ \|u_1\|_{H^{[n/2]}}+ \|u_1\|_{L^1}\le \epsilon, $$
we have a uniquely determined global (in time) small data Sobolev solution
$$ (u,v) \in \Big(\mathcal{C}\big([0,\ity),H^1 \cap L^\ity\big)\Big)^2 $$
to (\ref{pt1.4}). The following estimates hold for $k=0,1$:
\begin{align*}
\big\|\nabla^k u(t,\cdot)\big\|_{L^2} &\lesssim (1+t)^{-\frac{n}{4}-\frac{k}{2}+ \sigma(p_c,q_c)} \ell(t) \|(u_0,u_1),(v_0,v_1)\|_{\mathcal{A}}, \\
\|u(t,\cdot)\|_{L^\ity} &\lesssim (1+t)^{-\frac{n}{2}+ \sigma(p_c,q_c)} \ell(t) \|(u_0,u_1),(v_0,v_1)\|_{\mathcal{A}}, \\
\big\|\nabla^k v(t,\cdot)\big\|_{L^2} &\lesssim (1+t)^{-\frac{n}{4}-\frac{k}{2}} \|(u_0,u_1),(v_0,v_1)\|_{\mathcal{A}}, \\
\|v(t,\cdot)\|_{L^\ity} &\lesssim (1+t)^{-\frac{n}{2}} \|(u_0,u_1),(v_0,v_1)\|_{\mathcal{A}},
\end{align*}
where
\begin{equation}
\sigma(p_c,q_c):= \frac{q_c-p_c}{p_c q_c-1} \label{LossDecay.1}
\end{equation}
and the weighted function $\ell(t)$ is defined by
\begin{equation}
\ell(t):= \begin{cases}
1 &\text{ if (\ref{existenceCond.1.1.4}) holds}, \\
\left(\displaystyle\frac{\mu_1\big(c (1+t)^{-\varepsilon}\big)}{\mu_2\big(c (1+t)^{-\varepsilon}\big)}\right)^{\frac{1}{q_c+1}} &\text{ if (\ref{existenceCond.2.1.4}) holds}, \\
\end{cases} \label{LossDecay.2}
\end{equation}
with a sufficiently small constant $\varepsilon>0$.
\edl

\bnx
\fontshape{n}
\selectfont
We want to point out that the constant $\sigma(p_c,q_c)$ and the weighted function $\ell(t)$ appearing in the estimates for solutions in Theorem \ref{GlobalExistence1.4} represent some loss of decay in comparison with the corresponding estimates from the linear Cauchy problem.
\enx

\bvd
\fontshape{n}
\selectfont
We give some examples of the moduli of continuity $\mu_1$ and $\mu_2$ such that the assumptions (\ref{existenceCond.1.1.4}) and (\ref{existenceCond.2.1.4}) in Theorem \ref{GlobalExistence1.4} hold.
\begin{itemize}[leftmargin=*]
\item The assumption (\ref{existenceCond.1.1.4}) is fulfilled if we choose $\mu_1$ and $\mu_2$ by one of the following moduli of continuity:
\begin{enumerate}
\item $\mu(s)= s^\alpha$ with $\alpha \in (0,1]$;
\item $\mu(s)= \big(\log(1+s)\big)^\alpha$ with $\alpha \in (0,1)$;
\item $\mu(0)=0$ and $\mu(s)= \Big(\log\frac{1}{s}\Big)^{-\alpha}$ with $\alpha>1$;
\item $\mu(0)=0$ and $\mu(s)= \Big(\log\frac{1}{s}\Big)^{-1}\Big(\log\log\frac{1}{s}\Big)^{-1}\cdots \Big(\underbrace{\log\cdots\log}_{m \text{ times }\log}\frac{1}{s}\Big)^{-\alpha}$ with $m \in \N,\, \alpha>1$.
\end{enumerate}
\item The assumption (\ref{existenceCond.2.1.4}) is fulfilled if we choose $\mu_1$ and $\mu_2$ by one of the following pairs of moduli of continuity:
\begin{enumerate}
\item $\qquad$ $\mu_1(0)=0$ and $\mu_1(s)= \Big(\log\frac{1}{s}\Big)^{-\alpha_1}$ with $\alpha_1\le 1$, \\
\par $\qquad$ $\mu_2(0)=0$ and $\mu_2(s)= \Big(\log\frac{1}{s}\Big)^{-\alpha_2}$ with $\alpha_2>1$, \\
or $\quad$ $\mu_1(0)=0$ and $\mu_1(s)= \Big(\log\frac{1}{s}\Big)^{-\alpha_1}$ with $\alpha_1> 1$, \\
\par $\qquad$ $\mu_2(0)=0$ and $\mu_2(s)= \Big(\log\frac{1}{s}\Big)^{-\alpha_2}$ with $\alpha_2\le 1$, \\
provided that 
$$ \frac{q_c}{q_c+1}\alpha_1+ \frac{1}{q_c+1}\alpha_2>1; $$
\item $\quad$ $\mu_1(0)=0$ and $\mu_1(s)= \Big(\log\frac{1}{s}\Big)^{-1}\Big(\log\log\frac{1}{s}\Big)^{-1}\cdots \Big(\underbrace{\log\cdots\log}_{m \text{ times }\log}\frac{1}{s}\Big)^{-\alpha_1}$ with $m \in \N,\, \alpha_1\le 1$, \\
\par $\quad$ $\mu_2(0)=0$ and $\mu_2(s)= \Big(\log\frac{1}{s}\Big)^{-1}\Big(\log\log\frac{1}{s}\Big)^{-1}\cdots \Big(\underbrace{\log\cdots\log}_{m \text{ times }\log}\frac{1}{s}\Big)^{-\alpha_2}$ with $m \in \N,\, \alpha_2>1$, \\
or\, $\mu_1(0)=0$ and $\mu_1(s)= \Big(\log\frac{1}{s}\Big)^{-1}\Big(\log\log\frac{1}{s}\Big)^{-1}\cdots \Big(\underbrace{\log\cdots\log}_{m \text{ times }\log}\frac{1}{s}\Big)^{-\alpha_1}$ with $m \in \N,\, \alpha_1> 1$, \\
\par $\quad$ $\mu_2(0)=0$ and $\mu_2(s)= \Big(\log\frac{1}{s}\Big)^{-1}\Big(\log\log\frac{1}{s}\Big)^{-1}\cdots \Big(\underbrace{\log\cdots\log}_{m \text{ times }\log}\frac{1}{s}\Big)^{-\alpha_2}$ with $m \in \N,\, \alpha_2\le 1$, \\
provided that
$$ \frac{q_c}{q_c+1}\alpha_1+ \frac{1}{q_c+1}\alpha_2>1. $$
\end{enumerate}
Intuitively, from the two latter examples one can think of the modulus of continuity
$$ \mu_{1,2}:= \mu_{1,2}(s)= \big(\mu_1(s)\big)^{\frac{q_c}{q_c+1}}\big(\mu_2(s)\big)^{\frac{1}{q_c+1}} $$
as a ``middle'' modulus of continuity between $\mu_1$ and $\mu_2$ being subject to the following conditions:
$$ \int_0^{c}\frac{\mu_1(s)}{s}ds= \ity\quad \text{ and }\quad \int_0^{c}\frac{\mu_2(s)}{s}ds< \ity $$
or
$$ \int_0^{c}\frac{\mu_1(s)}{s}ds< \ity\quad \text{ and }\quad \int_0^{c}\frac{\mu_2(s)}{s}ds= \ity. $$
Then, we may take, among other things, suitable choices of $\mu_1$ as well as $\mu_2$ to claim that $\mu_{1,2}$ satisfies
$$ \int_0^{c}\frac{\mu_{1,2}(s)}{s}ds< \ity. $$
\end{itemize}
\evd

\bdl[\textbf{Blow-up}] \label{dloptimal1}
Assume that the initial data $u_0=v_0=0$ and $u_1,v_1 \in L^1$ satisfy the following relations:
\begin{equation} \label{optimaldata}
\int_{\R^n} u_1(x)\,dx >0 \quad \text{ and }\quad \int_{\R^n} v_1(x)\,dx >0.
\end{equation}
Moreover, we suppose the following assumptions of moduli of continuity:
\begin{equation}
s^k\mu_j^{(k)}(s)= o\big(\mu_j(s)\big) \quad \text{ as }s \to +0 \text{ with }j,k=1,2, \label{**}
\end{equation}
and
\begin{equation}
\int_0^{c}\frac{1}{s}\big(\mu_1(s)\big)^{\frac{q_c}{q_c+1}}\big(\mu_2(s)\big)^{\frac{1}{q_c+1}}\,ds= \ity, \label{blowupCon}
\end{equation}
where $c>0$ is a suitable small constant. Then, there is no global (in time) Sobolev solution to (\ref{pt1.4}).
\edl

\bvd
\fontshape{n}
\selectfont
We give some examples of the moduli of continuity $\mu_1$ and $\mu_2$ fulfilling the assumption (\ref{blowupCon}) in Theorem \ref{dloptimal1} as follows:
\begin{enumerate}[leftmargin=*]
\item $\qquad$ $\mu_1(0)=0$ and $\mu_1(s)= \Big(\log\frac{1}{s}\Big)^{-\alpha_1}$ with $\alpha_1> 0$, \\
\par $\quad$ $\mu_2(0)=0$ and $\mu_2(s)= \Big(\log\frac{1}{s}\Big)^{-\alpha_2}$ with $\alpha_2> 0$, \\
provided that 
$$ \frac{q_c}{q_c+1}\alpha_1+ \frac{1}{q_c+1}\alpha_2\le 1; $$
\item $\quad$ $\mu_1(0)=0$ and $\mu_1(s)= \Big(\log\frac{1}{s}\Big)^{-1}\Big(\log\log\frac{1}{s}\Big)^{-1}\cdots \Big(\underbrace{\log\cdots\log}_{m \text{ times }\log}\frac{1}{s}\Big)^{-\alpha_1}$ with $m \in \N,\, \alpha_1> 0$, \\
\par $\mu_2(0)=0$ and $\mu_2(s)= \Big(\log\frac{1}{s}\Big)^{-1}\Big(\log\log\frac{1}{s}\Big)^{-1}\cdots \Big(\underbrace{\log\cdots\log}_{m \text{ times }\log}\frac{1}{s}\Big)^{-\alpha_2}$ with $m \in \N,\, \alpha_2> 0$, \\
provided that
$$ \frac{q_c}{q_c+1}\alpha_1+ \frac{1}{q_c+1}\alpha_2\le 1. $$
\end{enumerate}
\evd

%===================================================Proofs
\section{Proofs of main results}
\subsection{Global existence}
In order to prove our global (in time) existence of small data Sobolev solutions, the following preliminary lemmas come into play.
\bbd[Lemma $1$ in \cite{Matsumura}] \label{bd1.4.1}
The solutions to the corresponding linear equations to (\ref{pt1.4}) satisfy the following estimates:
\begin{align*}
\big\|\nabla^k w(t,\cdot)\big\|_{L^2} &\lesssim (1+t)^{-\frac{n}{4}- \frac{k}{2}}\big(\|w_0\|_{L^1}+ \|w_0\|_{H^k}+ \|w_1\|_{L^1}+ \|w_1\|_{H^{k-1}}\big),
\end{align*}
with $k=0,1,1+[n/2]$ and
\begin{align*}
\|w(t,\cdot)\|_{L^\ity} &\lesssim (1+t)^{-\frac{n}{2}}\big(\|w_0\|_{L^1}+ \|w_0\|_{H^{1+[n/2]}}+ \|w_1\|_{L^1}+ \|w_1\|_{H^{[n/2]}}\big),
\end{align*}
where $w$ stands for $u$ or $v$.
\ebd

\bbd \label{bd1.4.2}
Let $\mu_1=\mu_1(s)$ and $\mu_2=\mu_2(s)$ are moduli of continuity. The following estimates hold:
\begin{align*}
\text{\fontshape{n}\selectfont(a)}&\,\, \int_0^t (1+t-\tau)^{-\alpha}(1+\tau)^{-1}\big(\mu_1\big(C(1+\tau)^{-\gamma}\big)\big)^{\beta_1}\big(\mu_2\big(C(1+\tau)^{-\gamma}\big)\big)^{\beta_2}d\tau \\
&\hspace{2cm} \lesssim (1+t)^{-\alpha}\int_0^t (1+\tau)^{-1}\big(\mu_1\big(C(1+\tau)^{-\gamma}\big)\big)^{\beta_1}\big(\mu_2\big(C(1+\tau)^{-\gamma}\big)\big)^{\beta_2}d\tau \\ 
&\text{ for any } \alpha\le 1 \text{ and for all } \beta_1,\beta_2,\gamma \ge 0, \\
\text{\fontshape{n}\selectfont(b)}&\,\, \int_0^\ity (1+\tau)^{-1}\big(\mu_1\big(C(1+\tau)^{-\gamma}\big)\big)^{\beta_1}\big(\mu_2\big(C(1+\tau)^{-\gamma}\big)\big)^{\beta_2}d\tau = C_0\int_0^{C} \frac{1}{s}\big(\mu_1(s)\big)^{\beta_1}\big(\mu_2(s)\big)^{\beta_2}ds \\
&\text{ for any } \beta_1,\beta_2 \ge 0 \text{ and for all } \gamma> 0, \text{ where } C_0 \text{ is a suitable positive constant.}
\end{align*}
\ebd

\begin{proof}
To prove the first estimate, we divide the left-hand side integral into two parts as follows:
\begin{align*}
& \int_0^t (1+t-\tau)^{-\alpha}(1+\tau)^{-1}\big(\mu_1\big(C(1+\tau)^{-\gamma}\big)\big)^{\beta_1}\big(\mu_2\big(C(1+\tau)^{-\gamma}\big)\big)^{\beta_2}d\tau \\ 
&\quad = \int_0^{t/2} (1+t-\tau)^{-\alpha}(1+\tau)^{-1}\big(\mu_1\big(C(1+\tau)^{-\gamma}\big)\big)^{\beta_1}\big(\mu_2\big(C(1+\tau)^{-\gamma}\big)\big)^{\beta_2}d\tau \\
&\qquad+ \int_{t/2}^t (1+t-\tau)^{-\alpha}(1+\tau)^{-1}\big(\mu_1\big(C(1+\tau)^{-\gamma}\big)\big)^{\beta_1}\big(\mu_2\big(C(1+\tau)^{-\gamma}\big)\big)^{\beta_2}d\tau \\
&\quad =: I_1+I_2.
\end{align*}
Using the relation $1+t-\tau \approx 1+t$ for any $\tau \in [0,t/2]$ one derives
\begin{align}
I_1 &\lesssim (1+t)^{-\alpha}\int_0^{t/2} (1+\tau)^{-1}\big(\mu_1\big(C(1+\tau)^{-\gamma}\big)\big)^{\beta_1}\big(\mu_2\big(C(1+\tau)^{-\gamma}\big)\big)^{\beta_2}d\tau \nonumber \\
&\le (1+t)^{-\alpha}\int_0^t (1+\tau)^{-1}\big(\mu_1\big(C(1+\tau)^{-\gamma}\big)\big)^{\beta_1}\big(\mu_2\big(C(1+\tau)^{-\gamma}\big)\big)^{\beta_2}d\tau. \label{lem1.4.2.1}
\end{align}
In addition, we notice the relation $1+\tau \approx 1+t \text{ for any }\tau \in [t/2,t]$ to deal with $I_2$ in the following way:
$$ I_2 \lesssim (1+t)^{-\alpha}\int_{t/2}^t (1+t-\tau)^{-\alpha}(1+\tau)^{\alpha-1}\big(\mu_1\big(C(1+\tau)^{-\gamma}\big)\big)^{\beta_1}\big(\mu_2\big(C(1+\tau)^{-\gamma}\big)\big)^{\beta_2}d\tau. $$
Due to the hypothesis $\alpha \le 1$ and $\gamma \ge 0$, it holds
$$ \begin{cases}
(1+\tau)^{\alpha-1} \le (1+t-\tau)^{\alpha-1} \\
(1+\tau)^{-\gamma} \le (1+t-\tau)^{-\gamma}
\end{cases} \quad \text{ for any } \tau \in [t/2,t]. $$
As a result, this leads to
\begin{align}
I_2 &\lesssim (1+t)^{-\alpha}\int_{t/2}^t (1+t-\tau)^{-1}\big(\mu_1\big(C(1+t-\tau)^{-\gamma}\big)\big)^{\beta_1}\big(\mu_2\big(C(1+t-\tau)^{-\gamma}\big)\big)^{\beta_2}d\tau \nonumber \\ 
&= (1+t)^{-\alpha}\int_0^{t/2} (1+\rho)^{-1}\big(\mu_1\big(C(1+\rho)^{-\gamma}\big)\big)^{\beta_1}\big(\mu_2\big(C(1+\rho)^{-\gamma}\big)\big)^{\beta_2}d\rho \nonumber \\
&\le (1+t)^{-\alpha}\int_0^t (1+\rho)^{-1}\big(\mu_1\big(C(1+\rho)^{-\gamma}\big)\big)^{\beta_1}\big(\mu_2\big(C(1+\rho)^{-\gamma}\big)\big)^{\beta_2}d\rho, \label{lem1.4.2.2}
\end{align}
where we have used the increasing property of the functions $\mu_1,\mu_2$ as well as change of variables $\rho=t-\tau$. Combining (\ref{lem1.4.2.1}) and (\ref{lem1.4.2.2}) we may conclude the estimate (a). Finally, a standard change of variables follows immediately the estimate (b).
\end{proof}

\bbd[Gagliardo-Nirenberg inequality, see \cite{ReissigEbert,Friedman}] \label{bd1.4.3}
Let $j,m \in \N$ with $j<m$ and $w \in \mathcal{C}^m_0(\R^n)$. Let us assume $\frac{j}{m}\le \theta\le 1$ and $1\le r,r_1,r_2\le \ity$ such that
$$ j-\frac{n}{r}= \Big(m-\frac{n}{r_1}\Big)\theta- \frac{n}{r_2}(1-\theta). $$
Then, it holds
$$ \big\|\nabla^j w\big\|_{L^r}\lesssim \big\|\nabla^m w\big\|^\theta_{L^{r_1}}\,\|w\|^{1-\theta}_{L^{r_2}}, $$
provided that $(m-\frac{n}{r_1})-j \notin \N$, that is, $\frac{n}{r_1}>m-j$ or $\frac{n}{r_1} \notin \N$.
\par If $(m-\frac{n}{r_1})-j \in \N$, then Gagliardo-Nirenberg inequality holds provided that $\frac{j}{m}\le \theta< 1$.
\ebd

\begin{proof}[Proof of Theorem \ref{GlobalExistence1.4}]
We introduce the solution space
$$ X(t):= \Big(\mathcal{C}\big([0,t],H^1 \cap L^\ity\big)\Big)^2 $$
with the norm
\begin{align*}
\|(u,v)\|_{X(t)}:= \sup_{0\le \tau \le t} \Big( (1+\tau)^{\frac{n}{4}- \sigma(p_c,q_c)} \big(\ell(\tau)\big)^{-1}\|u(\tau,\cdot)\|_{L^2} &+ (1+\tau)^{\frac{n}{4}+\frac{1}{2}- \sigma(p_c,q_c)} \big(\ell(\tau)\big)^{-1}\big\|\nabla u(\tau,\cdot)\big\|_{L^2} \\
&+ (1+\tau)^{\frac{n}{2}- \sigma(p_c,q_c)} \big(\ell(\tau)\big)^{-1}\|u(\tau,\cdot)\|_{L^\ity} \\
+ (1+\tau)^{\frac{n}{4}}\|v(\tau,\cdot)\|_{L^2} &+ (1+\tau)^{\frac{n}{4}+\frac{1}{2}}\big\|\nabla v(\tau,\cdot)\big\|_{L^2} \\
&+ (1+\tau)^{\frac{n}{2}}\|v(\tau,\cdot)\|_{L^\ity} \Big),
\end{align*}
where $\sigma(p_c,q_c)$ and the weighted function $\ell(\tau)$ are determined as in (\ref{LossDecay.1}) and (\ref{LossDecay.2}), respectively. We denote by $\mathcal{K}_0(t,x)$ and $\mathcal{K}_1(t,x)$ the fundamental solutions to the corresponding linear Cauchy problems for (\ref{pt1.4}). Then, the solutions to (\ref{pt1.4}) with vanishing right-hand side are written by
$$ \begin{cases}
u^{\text{ln}}(t,x)= \mathcal{K}_0(t,x) \ast_{x} u_0(x)+ \mathcal{K}_1(t,x) \ast_{x} u_1(x), \\
v^{\text{ln}}(t,x)= \mathcal{K}_0(t,x) \ast_{x} v_0(x)+ \mathcal{K}_1(t,x) \ast_{x} v_1(x).
\end{cases} $$
Thanks to Duhamel's principle, the formal integral representation of solutions to (\ref{pt1.4}) is
$$ \begin{cases}
u(t,x)= u^{\text{ln}}(t,x)+ \displaystyle\int_0^t \mathcal{K}_1(t-\tau,x) \ast_x \big(|v(\tau,x)|^{p_c}\mu_1(|v(\tau,x)|)\big) d\tau=: u^{\text{ln}}(t,x)+ u^{\text{nl}}(t,x), \\
\\
v(t,x)= v^{\text{ln}}(t,x)+ \displaystyle\int_0^t \mathcal{K}_1(t-\tau,x) \ast_x \big(|u(\tau,x)|^{q_c}\mu_2(|u(\tau,x)|)\big) d\tau=: v^{\text{ln}}(t,x)+ v^{\text{nl}}(t,x).
\end{cases} $$
For all $t>0$, we define the operator
$$ \Psi: \quad (u,v) \in X(t) \mapsto \Psi(u,v)(t,x)= \big(u^{\text{ln}}(t,x)+ u^{\text{nl}}(t,x), v^{\text{ln}}(t,x)+ v^{\text{nl}}(t,x)\big). $$
Our aim is use Banach's fixed point theorem to arrive at global (in time) existence of small data solutions to (\ref{pt1.4}). To establish this, we need to indicate that the operator $\Psi$ satisfies the following two inequalities:
\begin{align}
\|\Psi(u,v)\|_{X(t)} &\lesssim \|(u_0,u_1),(v_0,v_1)\|_{\mathcal{A}}+ \|(u,v)\|^{p_c}_{X(t)}+ \|(u,v)\|^{q_c}_{X(t)}, \label{pt3.3}\\
\|\Psi(u,v)-\Psi(\bar{u},\bar{v})\|_{X(t)} &\lesssim \|(u,v)-(\bar{u},\bar{v})\|_{X(t)} \Big(\|(u,v)\|^{p_c-1}_{X(t)}+ \|(\bar{u},\bar{v})\|^{p_c-1}_{X(t)} \nonumber \\
&\hspace{6cm}+ \|(u,v)\|^{q_c-1}_{X(t)}+ \|(\bar{u},\bar{v})\|^{q_c-1}_{X(t)}\Big). \label{pt3.4}
\end{align}
At first, it is clear to conclude that the estimate
$$ \big\|(u^{\text{ln}}, v^{\text{ln}})\big\|_{X(t)} \lesssim \|(u_0,u_1),(v_0,v_1)\|_{\mathcal{A}} $$
holds from Lemma \ref{bd1.4.1}. Hence, it suffices to prove only the following inequality instead of (\ref{pt3.3}):
\begin{equation}
\big\|(u^{\text{nl}}, v^{\text{nl}})\big\|_{X(t)} \lesssim \|(u,v)\|^{p_c}_{X(t)}+ \|(u,v)\|^{q_c}_{X(t)}. \label{pt3.31}
\end{equation}
Before indicating the inequality (\ref{pt3.31}), we require the following auxiliary estimates:
\begin{align}
&\big\||v(\tau,\cdot)|^{p_c}\mu_1\big(|v(\tau,\cdot)|\big)\big\|_{L^1 \cap L^2} \nonumber \\
&\qquad \lesssim (1+\tau)^{-1+\frac{q_c-p_c}{p_c q_c-1}}\,\,\mu_1\Big(c (1+ \tau)^{-\frac{1+q_c}{p_c q_c-1}}\Big)\|(u,v)\|^{p_c}_{X(t)}, \label{t3.1} \\
&\big\||u(\tau,\cdot)|^{q_c}\mu_2\big(|u(\tau,\cdot)|\big)\big\|_{L^1 \cap L^2} \nonumber \\
&\qquad \lesssim (1+\tau)^{-1} \big(\ell(\tau)\big)^{q_c}\,\,\mu_2\Big(c (1+ \tau)^{-\frac{1+p_c}{p_c q_c-1}} \ell(\tau)\Big)\|(u,v)\|^{q_c}_{X(t)}, \label{t3.2} \\
&\big\||v(\tau,\cdot)|^{p_c}\mu_1\big(|v(\tau,\cdot)|\big)\big\|_{L^1 \cap H^1} \nonumber \\
&\qquad \lesssim (1+\tau)^{-1+\frac{q_c-p_c}{p_c q_c-1}}\,\,\mu_1\Big(c (1+ \tau)^{-\frac{1+q_c}{p_c q_c-1}}\Big)\|(u,v)\|^{p_c}_{X(t)}, \label{t3.3} \\
&\big\||u(\tau,\cdot)|^{q_c}\mu_2\big(|u(\tau,\cdot)|\big)\big\|_{L^1 \cap H^1} \nonumber \\
&\qquad \lesssim (1+\tau)^{-1} \big(\ell(\tau)\big)^{q_c}\,\,\mu_2\Big(c (1+ \tau)^{-\frac{1+p_c}{p_c q_c-1}} \ell(\tau)\Big)\|(u,v)\|^{q_c}_{X(t)}, \label{t3.4}
\end{align}
where $c >0$ is a suitable small constant. Indeed, we may re-write
\begin{align*}
\big\||v(\tau,\cdot)|^{p_c}\big\|_{L^1 \cap L^2} &= \|v(\tau,\cdot)\|^{p_c}_{L^{p_c}}+ \|v(\tau,\cdot)\|^{p_c}_{L^{2p_c}}, \\ 
\big\||u(\tau,\cdot)|^{q_c}\big\|_{L^1 \cap L^2} &= \|u(\tau,\cdot)\|^{q_c}_{L^{q_c}}+ \|u(\tau,\cdot)\|^{q_c}_{L^{2q_c}}.
\end{align*}
The application of Gagliardo-Nirenberg inequality from Lemma \ref{bd1.4.3} to control the norms
$$ \|v(\tau,\cdot)\|_{L^{p_c}},\quad \|v(\tau,\cdot)\|_{L^{2p_c}},\quad \|u(\tau,\cdot)\|_{L^{q_c}} \quad \text{ and }\quad \|u(\tau,\cdot)\|^{q_c}_{L^{2q_c}} $$
follows immediately
\begin{align*}
\big\||v(\tau,\cdot)|^{p_c}\big\|_{L^1 \cap L^2} &\lesssim (1+\tau)^{-\frac{n}{2}(p_c- 1)}\|(u,v)\|^{p_c}_{X(t)},  \\ 
\big\||u(\tau,\cdot)|^{q_c}\big\|_{L^1 \cap L^2} &\lesssim (1+\tau)^{-\frac{n}{2}(q_c- 1)+ \sigma(p_c,q_c)q_c} \big(\ell(\tau)\big)^{q_c}\|(u,v)\|^{q_c}_{X(t)}.
\end{align*}
Since $\mu_1,\mu_2$ are increasing functions, from the definition of the norm in $X(t)$ one obtains
\begin{align*}
\big\|\mu_1\big(|v(\tau,\cdot)|\big)\big\|_{L^\ity}\le \mu_1\big(\|v(\tau,\cdot)\|_{L^\ity}\big) &\le \mu_1\Big(C (1+ \tau)^{-\frac{n}{2}} \|(u,v)\|_{X(t)}\Big) \\
&\le \mu_1\Big(c (1+ \tau)^{-\frac{n}{2}}\Big)
\end{align*}
and
\begin{align*}
\big\|\mu_2\big(|u(\tau,\cdot)|\big)\big\|_{L^\ity}\le \mu_2\big(\|u(\tau,\cdot)\|_{L^\ity}\big)&\le \mu_2\Big(C (1+ \tau)^{-\frac{n}{2}+ \sigma(p_c,q_c)} \ell(\tau) \|(u,v)\|_{X(t)}\Big) \\
&\le \mu_2\Big(c (1+ \tau)^{-\frac{n}{2}+ \sigma(p_c,q_c)} \ell(\tau)\Big)
\end{align*}
with $c:= C\e_0$, where $\e_0$ is a sufficiently small constant such that $\|(u,v)\|_{X(t)}\le \e_0$. For this reason, we may arrive at
\begin{align*}
\big\||v(\tau,\cdot)|^{p_c}\mu_1\big(|v(\tau,\cdot)|\big)\big\|_{L^1 \cap L^2} &\lesssim \big\||v(\tau,\cdot)|^{p_c}\big\|_{L^1 \cap L^2}\,\, \big\|\mu_1\big(|v(\tau,\cdot)|\big)\big\|_{L^\ity} \\
&\lesssim (1+\tau)^{-\frac{n}{2}(p_c- 1)}\,\,\mu_1\Big(c (1+ \tau)^{-\frac{n}{2}}\Big)\|(u,v)\|^{p_c}_{X(t)} \\
&= (1+\tau)^{-1+\frac{q_c-p_c}{p_c q_c-1}}\,\,\mu_1\Big(c (1+ \tau)^{-\frac{1+q_c}{p_c q_c-1}}\Big)\|(u,v)\|^{p_c}_{X(t)}
\end{align*}
and
\begin{align*}
&\big\||u(\tau,\cdot)|^{q_c}\mu_2\big(|u(\tau,\cdot)|\big)\big\|_{L^1 \cap L^2} \\
&\qquad \lesssim \big\||u(\tau,\cdot)|^{q_c}\big\|_{L^1 \cap L^2}\,\, \big\|\mu_2\big(|u(\tau,\cdot)|\big)\big\|_{L^\ity} \\
&\qquad \lesssim (1+\tau)^{-\frac{n}{2}(q_c- 1)+ \sigma(p_c,q_c)q_c} \big(\ell(\tau)\big)^{q_c}\,\,\mu_2\Big(c (1+ \tau)^{-\frac{n}{2}+ \sigma(p_c,q_c)} \ell(\tau)\Big)\|(u,v)\|^{q_c}_{X(t)} \\
&\qquad= (1+\tau)^{-1} \big(\ell(\tau)\big)^{q_c}\,\,\mu_2\Big(c (1+ \tau)^{-\frac{1+p_c}{p_c q_c-1}} \ell(\tau)\Big)\|(u,v)\|^{q_c}_{X(t)},
\end{align*}
where we notice that the relations
$$ -\frac{n}{2}(p_c- 1)= -1+\frac{q_c-p_c}{p_c q_c-1}, \quad -\frac{n}{2}(q_c- 1)+ \sigma(p_c,q_c)q_c=-1, \quad -\frac{n}{2}+ \sigma(p_c,q_c)= -\frac{1+p_c}{p_c q_c-1} $$
hold by (\ref{CriticalCurve}) and (\ref{LossDecay.1}). This completes the proof of (\ref{t3.1}) and (\ref{t3.2}). In order to show the two remaining estimates (\ref{t3.3}) and (\ref{t3.4}), we have the following expressions:
\begin{align}
\big\||v(\tau,\cdot)|^{p_c}\mu_1\big(|v(\tau,\cdot)|\big)\big\|_{L^1 \cap H^1} &= \big\||v(\tau,\cdot)|^{p_c}\mu_1\big(|v(\tau,\cdot)|\big)\big\|_{L^1 \cap L^2}+ \big\||v(\tau,\cdot)|^{p_c}\mu_1\big(|v(\tau,\cdot)|\big)\big\|_{\dot{H}^1}, \label{t3.5} \\ 
\big\||u(\tau,\cdot)|^{q_c}\mu_2\big(|u(\tau,\cdot)|\big)\big\|_{L^1 \cap H^1} &= \big\||u(\tau,\cdot)|^{q_c}\mu_2\big(|u(\tau,\cdot)|\big)\big\|_{L^1 \cap L^2}+ \big\||u(\tau,\cdot)|^{q_c}\mu_2\big(|u(\tau,\cdot)|\big)\big\|_{\dot{H}^1}. \label{t3.6} 
\end{align}
Therefore, it is reasonable to control the two additional norms only
$$ \big\||v(\tau,\cdot)|^{p_c}\mu_1\big(|v(\tau,\cdot)|\big)\big\|_{\dot{H}^1} \quad \text{ and }\quad \big\||u(\tau,\cdot)|^{q_c}\mu_2\big(|u(\tau,\cdot)|\big)\big\|_{\dot{H}^1}. $$
Observing that
\begin{align*}
\nabla \big(|v(\tau,x)|^{p_c}\mu_1(|v(\tau,x)|)\big) &= p_c |v(\tau,x)|^{p_c-2}\, v(\tau,x)\, \nabla v(\tau,x)\, \mu_1\big(|v(\tau,x)|\big) \\ 
&\quad + |v(\tau,x)|^{p_c}\, \mu_1'\big(|v(\tau,x)|\big)\, \nabla |v(\tau,x)|
\end{align*}
one derives the relation
$$ \big|\nabla \big(|v(\tau,x)|^{p_c}\mu_1(|v(\tau,x)|)\big)\big| \lesssim |v(\tau,x)|^{p_c-1}\, \mu_1\big(|v(\tau,x)|\big)\, |\nabla v(\tau,x)| $$
since the assumption (\ref{*}). Thus, it follows that
\begin{align}
\big\||v(\tau,\cdot)|^{p_c}\mu_1\big(|v(\tau,\cdot)|\big)\big\|_{\dot{H}^1} &\lesssim \|v(\tau,\cdot)\|^{p_c-1}_{L^\ity} \big\|\mu_1\big(|v(\tau,\cdot)|\big)\big\|_{L^\ity} \|\nabla v(\tau,\cdot)\|_{L^2} \nonumber \\ 
&\lesssim (1+\tau)^{-\frac{n}{2}(p_c- 1)-\frac{n}{4}-\frac{1}{2}} \mu_1\Big(c (1+ \tau)^{-\frac{n}{2}}\Big) \|(u,v)\|^{p_c}_{X(t)}. \label{t3.7} 
\end{align}
In the same way we obtain
\begin{align}
&\big\||u(\tau,\cdot)|^{q_c}\mu_2\big(|u(\tau,\cdot)|\big)\big\|_{\dot{H}^1} \nonumber \\ 
&\qquad \lesssim (1+\tau)^{-\frac{n}{2}(q_c- 1)-\frac{n}{4}-\frac{1}{2}+ \sigma(p_c,q_c)q_c} \big(\ell(\tau)\big)^{q_c}\,\,\mu_2\Big(c (1+ \tau)^{-\frac{n}{2}+ \sigma(p_c,q_c)} \ell(\tau)\Big)\|(u,v)\|^{q_c}_{X(t)}. \label{t3.8}
\end{align}
Collecting (\ref{t3.1}), (\ref{t3.5}), (\ref{t3.7}) and (\ref{t3.2}), (\ref{t3.6}), (\ref{t3.8}) we may conclude (\ref{t3.3}) and (\ref{t3.4}), respectively. \medskip

Let us come back to show the inequality (\ref{pt3.31}). Our strategy is to use the estimates from Lemma \ref{bd1.4.1} and the derived estimates from (\ref{t3.1}) to (\ref{t3.4}) to achieve the following estimates for $k=0,1$:
\begin{align*}
\big\|\nabla^k u^{\text{nl}}(\tau,\cdot)\big\|_{L^2}&\lesssim \int_0^t (1+t-\tau)^{-\frac{n}{4}-\frac{k}{2}}\big\||v(\tau,\cdot)|^{p_c}\mu_1\big(|v(\tau,\cdot)|\big)\big\|_{L^1 \cap L^2}d\tau \\
&\lesssim \|(u,v)\|^{p_c}_{X(t)}\int_0^t (1+t-\tau)^{-\frac{n}{4}-\frac{k}{2}}(1+\tau)^{-1+\frac{q_c-p_c}{p_c q_c-1}}\,\, \mu_1\Big(c (1+ \tau)^{-\frac{1+q_c}{p_c q_c-1}}\Big)d\tau, \\
\big\|u^{\text{nl}}(t,\cdot)\big\|_{L^\ity} &\lesssim 
\begin{cases}
\displaystyle\int_0^t (1+t-\tau)^{-\frac{n}{2}}\big\||v(\tau,\cdot)|^{p_c}\mu_1\big(|v(\tau,\cdot)|\big)\big\|_{L^1 \cap L^2}d\tau &\text{ if }n=1 \\
\\
\displaystyle\int_0^t (1+t-\tau)^{-\frac{n}{2}}\big\||v(\tau,\cdot)|^{p_c}\mu_1\big(|v(\tau,\cdot)|\big)\big\|_{L^1 \cap H^1}d\tau &\text{ if }n=2 \\
\end{cases}
\end{align*}
\begin{align*}
&\lesssim \|(u,v)\|^{p_c}_{X(t)}\int_0^t (1+t-\tau)^{-\frac{n}{2}}(1+\tau)^{-1+\frac{q_c-p_c}{p_c q_c-1}}\,\, \mu_1\Big(c (1+ \tau)^{-\frac{1+q_c}{p_c q_c-1}}\Big)d\tau,
\end{align*}
and
\begin{align*}
\big\|\nabla^k v^{\text{nl}}(\tau,\cdot)\big\|_{L^2} &\lesssim \int_0^t (1+t-\tau)^{-\frac{n}{4}-\frac{k}{2}}\big\||u(\tau,\cdot)|^{q_c}\mu_2\big(|u(\tau,\cdot)|\big)\big\|_{L^1 \cap L^2}d\tau \\ 
&\lesssim \|(u,v)\|^{q_c}_{X(t)}\int_0^t (1+t-\tau)^{-\frac{n}{4}-\frac{k}{2}}(1+\tau)^{-1} \big(\ell(\tau)\big)^{q_c}\mu_2\Big(c (1+ \tau)^{-\frac{1+p_c}{p_c q_c-1}} \ell(\tau)\Big)d\tau, \\
\big\|v^{\text{nl}}(t,\cdot)\big\|_{L^\ity} &\lesssim 
\begin{cases}
\displaystyle\int_0^t (1+t-\tau)^{-\frac{n}{2}}\big\||u(\tau,\cdot)|^{q_c}\mu_2\big(|u(\tau,\cdot)|\big)\big\|_{L^1 \cap L^2}d\tau &\text{ if }n=1 \\
\\
\displaystyle\int_0^t (1+t-\tau)^{-\frac{n}{2}}\big\||u(\tau,\cdot)|^{q_c}\mu_2\big(|u(\tau,\cdot)|\big)\big\|_{L^1 \cap H^1}d\tau &\text{ if }n=2 \\ 
\end{cases} \\
&\lesssim \|(u,v)\|^{q_c}_{X(t)}\int_0^t (1+t-\tau)^{-\frac{n}{2}}(1+\tau)^{-1} \big(\ell(\tau)\big)^{q_c}\mu_2\Big(c (1+ \tau)^{-\frac{1+p_c}{p_c q_c-1}} \ell(\tau)\Big)d\tau.
\end{align*}
Let us now divide our consideration into two cases according to (\ref{LossDecay.2}) as follows:
\begin{itemize}[leftmargin=*]
\item \textbf{Case 1:} If the assumption (\ref{existenceCond.1.1.4}) holds, then we take $\ell(\tau)\equiv 1$. For this reason, we can proceed as follows:
\begin{align*}
\big\|\nabla^k u^{\text{nl}}(\tau,\cdot)\big\|_{L^2} &\lesssim (1+t)^{\frac{q_c-p_c}{p_c q_c-1}}\|(u,v)\|^{p_c}_{X(t)}\int_0^t (1+t-\tau)^{-\frac{n}{4}-\frac{k}{2}}(1+\tau)^{-1}\,\, \mu_1\Big(c (1+ \tau)^{-\frac{1+q_c}{p_c q_c-1}}\Big)d\tau \\
&\lesssim (1+t)^{-\frac{n}{4}-\frac{k}{2}+\frac{q_c-p_c}{p_c q_c-1}}\|(u,v)\|^{p_c}_{X(t)}\int_0^t (1+\tau)^{-1}\,\, \mu_1\Big(c (1+ \tau)^{-\frac{1+q_c}{p_c q_c-1}}\Big)d\tau \\
&\le C (1+t)^{-\frac{n}{4}-\frac{k}{2}+\frac{q_c-p_c}{p_c q_c-1}}\|(u,v)\|^{p_c}_{X(t)}\int_0^\ity (1+\tau)^{-1}\,\, \mu_1\Big(c (1+ \tau)^{-\frac{1+q_c}{p_c q_c-1}}\Big)d\tau \\
&= C (1+t)^{-\frac{n}{4}-\frac{k}{2}+\frac{q_c-p_c}{p_c q_c-1}}\|(u,v)\|^{p_c}_{X(t)}\int_0^c \frac{\mu_1(s)}{s}ds \\
&\lesssim (1+t)^{-\frac{n}{4}-\frac{k}{2}+\frac{q_c-p_c}{p_c q_c-1}}\|(u,v)\|^{p_c}_{X(t)},
\end{align*}
where we have applied Lemma \ref{bd1.4.2} by choosing
$$ \alpha= \frac{n}{4}+\frac{k}{2},\quad \beta_1=1,\quad \beta_2=0, \quad \gamma= \frac{1+q_c}{p_c q_c-1} $$
and used the assumption (\ref{existenceCond.1.1.4}) as well. Moreover, one may estimate
\begin{align*}
\big\|\nabla^k v^{\text{nl}}(\tau,\cdot)\big\|_{L^2} &\lesssim \|(u,v)\|^{q_c}_{X(t)}\int_0^t (1+t-\tau)^{-\frac{n}{4}-\frac{k}{2}}(1+\tau)^{-1}\,\, \mu_2\Big(c (1+ \tau)^{-\frac{1+p_c}{p_c q_c-1}}\Big)d\tau \\
&\lesssim (1+t)^{-\frac{n}{4}-\frac{k}{2}}\|(u,v)\|^{q_c}_{X(t)}\int_0^t (1+\tau)^{-1}\,\, \mu_2\Big(c (1+ \tau)^{-\frac{1+p_c}{p_c q_c-1}}\Big)d\tau \\
&\le C (1+t)^{-\frac{n}{4}-\frac{k}{2}}\|(u,v)\|^{q_c}_{X(t)}\int_0^\ity (1+\tau)^{-1}\,\, \mu_2\Big(c (1+ \tau)^{-\frac{1+p_c}{p_c q_c-1}}\Big)d\tau \\
&= C (1+t)^{-\frac{n}{4}-\frac{k}{2}}\|(u,v)\|^{q_c}_{X(t)}\int_0^c \frac{\mu_2(s)}{s}ds \\
&\lesssim (1+t)^{-\frac{n}{4}-\frac{k}{2}}\|(u,v)\|^{q_c}_{X(t)},
\end{align*}
where we have applied Lemma \ref{bd1.4.2} by choosing 
$$ \alpha= \frac{n}{4}+\frac{k}{2},\quad \beta_1=0,\quad \beta_2=1,\quad \gamma= \frac{1+p_c}{p_c q_c-1} $$
and used the assumption (\ref{existenceCond.1.1.4}) as well. Analogously, we also obtain the following estimates:
\begin{align*}
\big\|u^{\text{nl}}(t,\cdot)\big\|_{L^\ity} &\lesssim (1+t)^{-\frac{n}{2}+\frac{q_c-p_c}{p_c q_c-1}}\|(u,v)\|^{p_c}_{X(t)}, \\ 
\big\|v^{\text{nl}}(t,\cdot)\big\|_{L^\ity} &\lesssim (1+t)^{-\frac{n}{2}}\|(u,v)\|^{q_c}_{X(t)},
\end{align*}
by the choice
$$ \alpha= \frac{n}{2}, \quad \beta_1=1, \quad \beta_2=0, \quad \gamma= \frac{1+q_c}{p_c q_c-1} $$
or
$$ \alpha= \frac{n}{2}, \quad \beta_1=0, \quad \beta_2=1, \quad \gamma= \frac{1+p_c}{p_c q_c-1} $$
according to these estimates for $u^{\text{nl}}$ or $v^{\text{nl}}$. From the definition of the norm in $X(t)$, collecting all the above estimates completes the inequality (\ref{pt3.31}).
\item \textbf{Case 2:} If the assumption (\ref{existenceCond.2.1.4}) holds, then we take
$$ \ell(\tau)= \left(\frac{\mu_1\big(c (1+ \tau)^{-\varepsilon}\big)}{\mu_2\big(c (1+ \tau)^{-\varepsilon}\big)}\right)^{\frac{1}{q_c+1}}. $$
Following some arguments as we did in Case $1$ we may estimate
\begin{align*}
\big\|\nabla^k u^{\text{nl}}(\tau,\cdot)\big\|_{L^2} &\lesssim \|(u,v)\|^{p_c}_{X(t)}\int_0^t (1+t-\tau)^{-\frac{n}{4}-\frac{k}{2}}(1+\tau)^{-1+\frac{q_c-p_c}{p_c q_c-1}}\ell(\tau) \big(\ell(\tau)\big)^{-1} \mu_1\Big(c (1+ \tau)^{-\varepsilon}\Big)d\tau \\
&\hspace{5cm} \big(\text{since } \mu_1 \text{ is an increasing function}\big) \\
&\lesssim (1+t)^{\frac{q_c-p_c}{p_c q_c-1}}\ell(t) \|(u,v)\|^{p_c}_{X(t)} \\
&\quad \times \int_0^t (1+t-\tau)^{-\frac{n}{4}-\frac{k}{2}}(1+\tau)^{-1}\,\,\Big(\mu_1\big(c (1+ \tau)^{-\varepsilon}\big)\Big)^{\frac{q_c}{q_c+1}}\Big(\mu_2\big(c (1+ \tau)^{-\varepsilon}\big)\Big)^{\frac{1}{q_c+1}}d\tau \\
&\hspace{5cm} \big(\text{by (\ref{WeightedFunct.1}) in Remark \ref{Remark.WeightedFunct}}\big) \\
&\lesssim (1+t)^{-\frac{n}{4}-\frac{k}{2}+\frac{q_c-p_c}{p_c q_c-1}}\ell(t)\|(u,v)\|^{p_c}_{X(t)} \\
&\quad \times \int_0^t (1+\tau)^{-1}\,\, \Big(\mu_1\big(c (1+ \tau)^{-\varepsilon}\big)\Big)^{\frac{q_c}{q_c+1}}\Big(\mu_2\big(c (1+ \tau)^{-\varepsilon}\big)\Big)^{\frac{1}{q_c+1}}d\tau \\
&\le C (1+t)^{-\frac{n}{4}-\frac{k}{2}+\frac{q_c-p_c}{p_c q_c-1}}\ell(t)\|(u,v)\|^{p_c}_{X(t)} \\
&\quad \times \int_0^\ity (1+\tau)^{-1}\,\, \Big(\mu_1\big(c (1+ \tau)^{-\varepsilon}\big)\Big)^{\frac{q_c}{q_c+1}}\Big(\mu_2\big(c (1+ \tau)^{-\varepsilon}\big)\Big)^{\frac{1}{q_c+1}}d\tau \\
&= C (1+t)^{-\frac{n}{4}-\frac{k}{2}+\frac{q_c-p_c}{p_c q_c-1}}\ell(t)\|(u,v)\|^{p_c}_{X(t)}\int_0^c \frac{1}{s}\big(\mu_1(s)\big)^{\frac{q_c}{q_c+1}}\big(\mu_2(s)\big)^{\frac{1}{q_c+1}}ds \\
&\lesssim (1+t)^{-\frac{n}{4}-\frac{k}{2}+\frac{q_c-p_c}{p_c q_c-1}}\ell(t)\|(u,v)\|^{p_c}_{X(t)},
\end{align*}
where we have applied Lemma \ref{bd1.4.2} by choosing
$$ \alpha= \frac{n}{4}+\frac{k}{2},\quad \beta_1=\frac{q_c}{q_c+1},\quad \beta_2=\frac{1}{q_c+1}, \quad \gamma= \varepsilon $$
as well as used the assumption (\ref{existenceCond.2.1.4}). In such a way one also has
\begin{align*}
\big\|\nabla^k v^{\text{nl}}(\tau,\cdot)\big\|_{L^2} &\lesssim \|(u,v)\|^{q_c}_{X(t)}\int_0^t (1+t-\tau)^{-\frac{n}{4}-\frac{k}{2}}(1+\tau)^{-1}\big(\ell(\tau)\big)^{q_c}\,\, \mu_2\big(c (1+ \tau)^{-\varepsilon}\big)d\tau \\
&\hspace{5cm} \big(\text{by (\ref{WeightedFunct.2}) in Remark \ref{Remark.WeightedFunct}}\big) \\
&= \|(u,v)\|^{q_c}_{X(t)}\int_0^t (1+t-\tau)^{-\frac{n}{4}-\frac{k}{2}}(1+\tau)^{-1} \\
&\hspace{3cm} \times \Big(\mu_1\big(c (1+ \tau)^{-\varepsilon}\big)\Big)^{\frac{q_c}{q_c+1}}\Big(\mu_2\big(c (1+ \tau)^{-\varepsilon}\big)\Big)^{\frac{1}{q_c+1}}d\tau \\
&\lesssim (1+t)^{-\frac{n}{4}-\frac{k}{2}}\|(u,v)\|^{q_c}_{X(t)} \\
&\qquad \times \int_0^t (1+\tau)^{-1}\,\, \Big(\mu_1\big(c (1+ \tau)^{-\varepsilon}\big)\Big)^{\frac{q_c}{q_c+1}}\Big(\mu_2\big(c (1+ \tau)^{-\varepsilon}\big)\Big)^{\frac{1}{q_c+1}}d\tau \\
&\le C (1+t)^{-\frac{n}{4}-\frac{k}{2}}\|(u,v)\|^{q_c}_{X(t)} \\
&\qquad \times \int_0^\ity (1+\tau)^{-1}\,\, \Big(\mu_1\big(c (1+ \tau)^{-\varepsilon}\big)\Big)^{\frac{q_c}{q_c+1}}\Big(\mu_2\big(c (1+ \tau)^{-\varepsilon}\big)\Big)^{\frac{1}{q_c+1}}d\tau \\
&= C (1+t)^{-\frac{n}{4}-\frac{k}{2}}\|(u,v)\|^{q_c}_{X(t)}\int_0^c \frac{1}{s}\big(\mu_1(s)\big)^{\frac{q_c}{q_c+1}}\big(\mu_2(s)\big)^{\frac{1}{q_c+1}}ds \\
&\lesssim (1+t)^{-\frac{n}{4}-\frac{k}{2}}\|(u,v)\|^{q_c}_{X(t)},
\end{align*}
where we have employed Lemma \ref{bd1.4.2} by choosing 
$$ \alpha= \frac{n}{4}+\frac{k}{2},\quad \beta_1=\frac{q_c}{q_c+1},\quad \beta_2=\frac{1}{q_c+1}, \quad \gamma= \varepsilon $$
and used the assumption (\ref{existenceCond.2.1.4}) as well. Similarly, we may derive the following estimates:
\begin{align*}
\big\|u^{\text{nl}}(t,\cdot)\big\|_{L^\ity} &\lesssim (1+t)^{-\frac{n}{2}+\frac{q_c-p_c}{p_c q_c-1}}\ell(t)\|(u,v)\|^{p_c}_{X(t)}, \\ 
\big\|v^{\text{nl}}(t,\cdot)\big\|_{L^\ity} &\lesssim (1+t)^{-\frac{n}{2}}\|(u,v)\|^{q_c}_{X(t)},
\end{align*}
by the choice
$$ \alpha= \frac{n}{2},\quad \beta_1=\frac{q_c}{q_c+1},\quad \beta_2=\frac{1}{q_c+1}, \quad \gamma= \varepsilon $$
according to these estimates for both $u^{\text{nl}}$ and $v^{\text{nl}}$. From the definition of the norm in $X(t)$, we combine all the above estimates to complete the inequality (\ref{pt3.31}).
\end{itemize}

Next, let us prove the inequality (\ref{pt3.4}). For two elements $(u,v)$ and $(\bar{u},\bar{v})$ from $X(t)$, it is obvious that
$$ \Psi(u,v)(t,x)- \Psi(\bar{u},\bar{v})(t,x)= \big(u^{\text{nl}}(t,x)- \bar{u}^{\text{nl}}(t,x), v^{\text{nl}}(t,x)- \bar{v}^{\text{nl}}(t,x)\big). $$
Then, we use the same strategies as in the proof of the inequality (\ref{pt3.31}) to gain the following estimates with $k= 0,1$:
\begin{align*}
&\big\|\nabla^k \big(u^{\text{nl}}- \bar{u}^{\text{nl}}\big)(\tau,\cdot)\big\|_{L^2} \\
&\qquad \lesssim \int_0^t (1+t-\tau)^{-\frac{n}{4}-\frac{k}{2}}\big\||v(\tau,\cdot)|^{p_c}\mu_1\big(|v(\tau,\cdot)|\big)- |\bar{v}(\tau,\cdot)|^{p_c}\mu_1\big(|\bar{v}(\tau,\cdot)|\big)\big\|_{L^1 \cap L^2}d\tau,
\end{align*}
\begin{align*}
&\big\|\big(u^{\text{nl}}- \bar{u}^{\text{nl}}\big)(t,\cdot)\big\|_{L^\ity} \\
&\qquad \lesssim 
\begin{cases}
\displaystyle\int_0^t (1+t-\tau)^{-\frac{n}{2}}\big\||v(\tau,\cdot)|^{p_c}\mu_1\big(|v(\tau,\cdot)|\big)- |\bar{v}(\tau,\cdot)|^{p_c}\mu_1\big(|\bar{v}(\tau,\cdot)|\big)\big\|_{L^1 \cap L^2}d\tau &\text{ if }n=1, \\
\\
\displaystyle\int_0^t (1+t-\tau)^{-\frac{n}{2}}\big\||v(\tau,\cdot)|^{p_c}\mu_1\big(|v(\tau,\cdot)|\big)- |\bar{v}(\tau,\cdot)|^{p_c}\mu_1\big(|\bar{v}(\tau,\cdot)|\big)\big\|_{L^1 \cap H^1}d\tau &\text{ if }n=2, \\
\end{cases}
\end{align*}
and
\begin{align*}
&\big\|\nabla^k \big(v^{\text{nl}}- \bar{v}^{\text{nl}}\big)(\tau,\cdot)\big\|_{L^2} \\
&\qquad \lesssim \int_0^t (1+t-\tau)^{-\frac{n}{4}-\frac{k}{2}}\big\||u(\tau,\cdot)|^{q_c}\mu_2\big(|u(\tau,\cdot)|\big)- |\bar{u}(\tau,\cdot)|^{q_c}\mu_2\big(|\bar{u}(\tau,\cdot)|\big)\big\|_{L^1 \cap L^2}d\tau, \\ 
&\big\|\big(v^{\text{nl}}- \bar{v}^{\text{nl}}\big)(t,\cdot)\big\|_{L^\ity} \\
&\qquad \lesssim 
\begin{cases}
\displaystyle\int_0^t (1+t-\tau)^{-\frac{n}{2}}\big\||u(\tau,\cdot)|^{q_c}\mu_2\big(|u(\tau,\cdot)|\big)- |\bar{u}(\tau,\cdot)|^{q_c}\mu_2\big(|\bar{u}(\tau,\cdot)|\big)\big\|_{L^1 \cap L^2}d\tau &\text{ if }n=1, \\
\\
\displaystyle\int_0^t (1+t-\tau)^{-\frac{n}{2}}\big\||u(\tau,\cdot)|^{q_c}\mu_2\big(|u(\tau,\cdot)|\big)- |\bar{u}(\tau,\cdot)|^{q_c}\mu_2\big(|\bar{u}(\tau,\cdot)|\big)\big\|_{L^1 \cap H^1}d\tau &\text{ if }n=2. \\
\end{cases}
\end{align*}
Applying the mean value theorem gives the following integral representation:
\begin{align*}
&|v(\tau,x)|^{p_c}\mu_1\big(|v(\tau,x)|\big)- |\bar{v}(\tau,x)|^{p_c}\mu_1\big(|\bar{v}(\tau,x)|\big)\\ 
&\qquad= \big(v(\tau,x)- \bar{v}(\tau,x)\big)\int_0^1 d_{|v|}G\big(\omega v(\tau,x)+(1-\omega)\bar{v}(\tau,x)\big)\,d\omega,
\end{align*}
where $G(v)= |v|^{p_c}\mu_1(|v|)$. Since the condition (\ref{*}) of moduli of continuity holds, one gets
$$ d_{|v|}G(v)= p_c |v|^{p_c- 1}\mu_1(|v|)+ |v|^{p_c}d_{|v|}\mu_1(|v|) \lesssim |v|^{p_c- 1}\mu_1(|v|). $$
Thus, it follows that
\begin{align*}
&|v(\tau,x)|^{p_c}\mu_1\big(|v(\tau,x)|\big)- |\bar{v}(\tau,x)|^{p_c}\mu_1\big(|\bar{v}(\tau,x)|\big)\\ 
&\qquad \lesssim \big(v(\tau,x)- \bar{v}(\tau,x)\big)\int_0^1 \big|\omega v(\tau,x)+(1-\omega)\bar{v}(\tau,x)\big|^{p_c- 1}\mu_1\big(\big|\omega v(\tau,x)+(1-\omega)\bar{v}(\tau,x)\big|\big)\,d\omega.
\end{align*}
Similarly, we also obtain
\begin{align*}
&|u(\tau,x)|^{q_c}\mu_2\big(|u(\tau,x)|\big)- |\bar{u}(\tau,x)|^{q_c}\mu_2\big(|\bar{u}(\tau,x)|\big)\\ 
&\qquad \lesssim \big(u(\tau,x)- \bar{u}(\tau,x)\big)\int_0^1 \big|\omega u(\tau,x)+(1-\omega)\bar{u}(\tau,x)\big|^{q_c- 1}\mu_2\big(\big|\omega u(\tau,x)+(1-\omega)\bar{u}(\tau,x)\big|\big)\,d\omega.
\end{align*}
By the aid of H\"{o}lder's inequality and following the same manner as in the proof of the inequality (\ref{pt3.31}), we may arrive at the inequality (\ref{pt3.4}). Summarizing, the proof of Theorem \ref{GlobalExistence1.4} is completed.
\end{proof}

\bnx \label{Remark.WeightedFunct}
\fontshape{n}
\selectfont
Here we want to underline that in the proof of Theorem \ref{GlobalExistence1.4} we have used the following auxiliary properties of the weighted function $\ell(\tau)$ in Case $2$:
\begin{align}
\text{\fontshape{n}\selectfont i)}&\,\, (1+\tau)^{\frac{q_c-p_c}{p_c q_c-1}}\ell(\tau) \quad \text{ is increasing}; \label{WeightedFunct.1} \\
\text{\fontshape{n}\selectfont ii)}&\,\, (1+ \tau)^{-\frac{1+p_c}{p_c q_c-1}}\ell(\tau) \le (1+ \tau)^{-\varepsilon}. \label{WeightedFunct.2}
\end{align}
Indeed, by change of variables $s=c (1+ \tau)^{-\varepsilon}$ we may re-write
\begin{align*}
f(\tau):= (1+\tau)^{\frac{q_c-p_c}{p_c q_c-1}}\ell(\tau) &= (1+\tau)^{\frac{q_c-p_c}{p_c q_c-1}}\left(\frac{\mu_1\big(c (1+ \tau)^{-\varepsilon}\big)}{\mu_2\big(c (1+ \tau)^{-\varepsilon}\big)}\right)^{\frac{1}{q_c+1}} \\ 
&= C\,s^{-\frac{q_c-p_c}{\varepsilon(p_c q_c-1)}}\left(\frac{\mu_1(s)}{\mu_2(s)}\right)^{\frac{1}{q_c+1}}= C\,\left(s^{-\frac{(q_c-p_c)(q_c+1)}{\varepsilon(p_c q_c-1)}}\frac{\mu_1(s)}{\mu_2(s)}\right)^{\frac{1}{q_c+1}}.
\end{align*}
For this reason, in order to prove that $f(\tau)$ is an increasing function, it suffices to verify that
$$ h_1(s):= s^{-\frac{(q_c-p_c)(q_c+1)}{\varepsilon(p_c q_c-1)}}\mu_1(s) $$
is a decreasing function due to the increasing property of the function $\mu_2$. One has
\begin{align*}
h'_1(s) &= -\frac{(q_c-p_c)(q_c+1)}{\varepsilon(p_c q_c-1)}s^{-\frac{(q_c-p_c)(q_c+1)}{\varepsilon(p_c q_c-1)}-1}\mu_1(s)+ s^{-\frac{(q_c-p_c)(q_c+1)}{\varepsilon(p_c q_c-1)}}\mu'_1(s) \\ 
&\le s^{-\frac{(q_c-p_c)(q_c+1)}{\varepsilon(p_c q_c-1)}-1}\mu_1(s) \left(-\frac{(q_c-p_c)(q_c+1)}{\varepsilon(p_c q_c-1)}+ C\right) \quad \big(\text{since } s\mu_1'(s) \le C\mu_1(s) \text{ from } (\ref{*})\big) \\
&\le 0,
\end{align*}
because of the choice of a sufficiently small constant $\varepsilon>0$. This is to claim the first statement (\ref{WeightedFunct.1}). In an analogous way, we may conclude the second one (\ref{WeightedFunct.2}).
\enx

%=============================================================
\subsection{Blow-up result}
To show this result, the following generalized Jensen's inequality comes into play.
\bbd[Lemma $8$ in \cite{EbertGirardiReissig}] \label{Jenseninequality}
Let $\eta= \eta(x)$ be a defined and nonnegative function almost everywhere on $\Omega$, provided that $\eta$ is positive in a set of positive measure. Then, for each convex function $h$ on $\R$ the following inequality holds:
$$ h\left(\frac{\displaystyle\int_\Omega f(x)\eta(x)\,dx}{\displaystyle\int_\Omega \eta(x)\,dx}\right) \le \frac{\displaystyle\int_\Omega h\big(f(x)\big)\eta(x)\,dx}{\displaystyle\int_\Omega \eta(x)\,dx}, $$
where $f$ is any nonnegative function satisfying all the above integrals are meaningful.
\ebd

\begin{proof}[Proof of Theorem \ref{dloptimal1}]
Our proof relies on the ideas from the recent paper of the second author and his collaborators, where the authors have devoted to the study of the single semi-linear damped wave equation (\ref{pt1.2}). First of all, we introduce a test function $\varphi= \varphi(\rho)$ fulfilling
$$ \varphi \in \mathcal{C}_0^\ity\big([0,\ity)\big) \text{ and }
\varphi (\rho)=\begin{cases}
1 &\quad \text{ if }\quad 0 \le \rho \le 1/2, \\
\text{decreasing } &\quad \text{ if }\quad 1/2 \le \rho \le 1, \\
0 &\quad \text{ if }\quad \rho \ge 1.
\end{cases} $$
Also, we introduce the function $\varphi^*= \varphi^*(\rho)$ as follows:
$$ \varphi^*(\rho)= \begin{cases}
0 &\quad \text{ if }\quad 0 \le \rho< 1/2, \\
\varphi(\rho) &\quad \text{ if }\quad 1/2 \le \rho< \ity.
\end{cases} $$
Let $R$ be a large parameter in $[0,\ity)$. We denote the two functions
$$ \phi_R(t,x)= \Big(\varphi\Big(\frac{t^2+|x|^4}{R^4}\Big)\Big)^{\nu+2} \quad \text{ and }\quad \phi^*_R(t,x)= \Big(\varphi^*\Big(\frac{t^2+|x|^4}{R^4}\Big)\Big)^{\nu+2}, $$
where the parameter $\nu>0$ will be fixed later. Then, we may observe that
\begin{align*}
&\text{supp} \,\phi_R \subset Q_R:= \big\{(t,x) \,\,:\,\, (t,|x|) \in \big[0,R^2\big] \times \big[0,R\big] \big\}, \\
&\text{supp} \,\phi^*_R \subset Q^*_R:= Q_R \,\setminus \,\big\{(t,x) \,\,:\,\, (t,|x|) \in \big[0,R^2/\sqrt{2}\big] \times \big[0,(R/\sqrt[4]{2}\big] \big\}.
\end{align*}
Now we define the following two functionals:
\begin{align*}
I_R:= &\int_0^{\ity}\int_{\R^n}|v(t,x)|^{p_c}\mu_1\big(|v(t,x)|\big) \phi_R(t,x)\,dxdt= \int_{Q_R}|v(t,x)|^{p_c}\mu_1\big(|v(t,x)|\big) \phi_R(t,x)\,d(x,t), \\
J_R:= &\int_0^{\ity}\int_{\R^n}|u(t,x)|^{q_c}\mu_2\big(|u(t,x)|\big) \phi_R(t,x)\,dxdt= \int_{Q_R}|u(t,x)|^{q_c}\mu_2\big(|u(t,x)|\big) \phi_R(t,x)\,d(x,t).
\end{align*}
Let us assume that $(u,v)= (u(t,x),v(t,x))$ is a global (in time) Sobolev solution to (\ref{pt1.4}). We multiply the left-hand sides of (\ref{pt1.4}) by $\phi_R=\phi_R(t,x)$ and integrate by parts to achieve
\begin{align}
0\le I_R &= -\int_{\R^n} u_1(x)\phi_R(0,x)\,dx+ \int_{Q_R}u(t,x) \big(\partial_t^2 \phi_R(t,x)- \Delta \phi_R(t,x)- \partial_t \phi_R(t,x) \big)\,d(x,t) \nonumber \\
&=: -\int_{\R^n} u_1(x)\phi_R(0,x)\,dx + I^*_R \label{t4.1}
\end{align}
and
\begin{align}
0\le J_R &= -\int_{\R^n} v_1(x)\phi_R(0,x)\,dx+ \int_{Q_R}v(t,x) \big(\partial_t^2 \phi_R(t,x)- \Delta \phi_R(t,x)- \partial_t \phi_R(t,x) \big)\,d(x,t) \nonumber \\
&=: -\int_{\R^n} v_1(x)\phi_R(0,x)\,dx + J^*_R. \label{t4.2}
\end{align}
To estimate $I^*_R$ and $J^*_R$, a straightforward calculation gives the following estimates:
\begin{align*}
\big|\partial_t \phi_R(t,x)\big| &\lesssim \frac{1}{R^2}\Big(\varphi^*\Big(\frac{t^2+|x|^4}{R^4}\Big)\Big)^{\nu+1}, \\
\big|\partial_t^2 \phi_R(t,x)\big| &\lesssim \frac{1}{R^4}\Big(\varphi^*\Big(\frac{t^2+|x|^4}{R^4}\Big)\Big)^{\nu}, \\
\big| \Delta \phi_R(t,x)\big| &\lesssim \frac{1}{R^2} \Big(\varphi^*\Big(\frac{t^2+|x|^4}{R^4}\Big)\Big)^{\nu},
\end{align*}
where we have used the support conditions of $\phi_R$ and $\phi^*_R$. As a consequence, we can proceed $I^*_R$ and $J^*_R$ as follows:
\begin{equation}
|I^*_R| \lesssim \frac{1}{R^2} \int_{Q_R}|u(t,x)|\, \Big(\varphi^*\Big(\frac{t^2+|x|^4}{R^4}\Big)\Big)^{\nu}\,d(x,t)= \frac{1}{R^2} \int_{Q_R}|u(t,x)|\, \big(\phi^*_R(t,x)\big)^{\frac{\nu}{\nu+2}}\,d(x,t) \label{t4.3}
\end{equation}
and
\begin{equation}
|J^*_R| \lesssim \frac{1}{R^2} \int_{Q_R}|v(t,x)|\, \Big(\varphi^*\Big(\frac{t^2+|x|^4}{R^4}\Big)\Big)^{\nu}\,d(x,t)= \frac{1}{R^2} \int_{Q_R}|v(t,x)|\, \big(\phi^*_R(t,x)\big)^{\frac{\nu}{\nu+2}}\,d(x,t). \label{t4.4}
\end{equation}
Let us now devote to the estimation for the above integrals. For this purpose, we define the two functions $\Phi_p(s)= s^{p_c}\mu_1(s)$ and $\Phi_q(s)= s^{q_c}\mu_2(s)$. We have
\begin{align}
&\Phi_q\Big( |u(t,x)|\, \big(\phi^*_R(t,x)\big)^{\frac{\nu}{\nu+2}}\Big) \nonumber \\
&\qquad = |u(t,x)|^{q_c}\, \big(\phi^*_R(t,x)\big)^{\frac{\nu q_c}{\nu+2}} \mu_2\Big( |u(t,x)|\, \big(\phi^*_R(t,x)\big)^{\frac{\nu}{\nu+2}}\Big) \nonumber \\
&\qquad \le |u(t,x)|^{q_c}\, \big(\phi^*_R(t,x)\big)^{\frac{\nu q_c}{\nu+2}} \mu_2\big( |u(t,x)|\big)= \Phi_q\big( |u(t,x)|\big)\, \big(\phi^*_R(t,x)\big)^{\frac{\nu q_c}{\nu+2}} \label{t4.5}
\end{align}
since $\mu=\mu(s)$ is an increasing function and it holds
$$0\le \big(\phi^*_R(t,x)\big)^{\frac{\nu}{\nu+2}} \le 1 $$
for any $\nu>0$. It is obvious from the assumption (\ref{**}) that
$$ \Phi_q''(s)= s^{q_c-2}\Big(q_c (q_c-1)\mu_2(s)+2q_c\, s\mu'_2(s)+ s^2 \mu''_2(s)\Big) \ge 0, $$
that is, $\Phi_q$ is a convex function on a small interval $(0,c_0]$ with a sufficiently small constant $c_0>0$. Additionally, we can choose a convex continuation of $\Phi_q$ outside this interval to guarantee that $\Phi_q$ is convex on $[0,\ity)$. The application of the generalized Jensen's inequality from Lemma \ref{Jenseninequality} with $h(s)= \Phi_q(s)$, $f(t,x)= |u(t,x)|\big(\phi^*_R(t,x)\big)^{\frac{\nu}{\nu+2}}$, $\eta \equiv 1$ and $\Omega \equiv Q^*_R$ leads to the following estimate:
\begin{align*}
&\Phi_q\left(\frac{\displaystyle\int_{Q^*_R} |u(t,x)|\big(\phi^*_R(t,x)\big)^{\frac{\nu}{\nu+2}}\,d(x,t)}{\displaystyle\int_{Q^*_R} 1\,d(x,t)}\right) \le \frac{\displaystyle\int_{Q^*_R} \Phi_q\Big( |u(t,x)|\big(\phi^*_R(t,x)\big)^{\frac{\nu}{\nu+2}}\Big)\,d(x,t)}{\displaystyle\int_{Q^*_R} 1\,d(x,t)}.
\end{align*}
It is clear to verify that
$$ \int_{Q^*_R} 1\,d(x,t)\approx R^{n+2}. $$
Thus, it follows
\begin{align}
\Phi_q\left(\frac{\displaystyle\int_{Q^*_R} |u(t,x)|\big(\phi^*_R(t,x)\big)^{\frac{\nu}{\nu+2}}\,d(x,t)}{R^{n+2}}\right) &\le \frac{\displaystyle\int_{Q^*_R} \Phi_q\Big( |u(t,x)|\big(\phi^*_R(t,x)\big)^{\frac{\nu}{\nu+2}}\Big)\,d(x,t)}{R^{n+2}} \nonumber \\
&\le \frac{\displaystyle\int_{Q_R} \Phi_q\Big( |u(t,x)|\big(\phi^*_R(t,x)\big)^{\frac{\nu}{\nu+2}}\Big)\,d(x,t)}{R^{n+2}}. \label{t4.6}
\end{align}
From the estimates (\ref{t4.5}) and (\ref{t4.6}) one derives
\begin{align}
&\Phi_q\left(\frac{\displaystyle\int_{Q^*_R} |u(t,x)|\big(\phi^*_R(t,x)\big)^{\frac{\nu}{\nu+2}}\,d(x,t)}{R^{n+2}}\right) \le \frac{\displaystyle\int_{Q_R} \Phi_q\big( |u(t,x)|\big)\, \big(\phi^*_R(t,x)\big)^{\frac{\nu q_c}{\nu+2}}\,d(x,t)}{R^{n+2}}. \label{t4.7}
\end{align}
Due to the fact that $\mu=\mu(s)$ is a strictly increasing function, $\Phi_q=\Phi_q(s)$ is also a strictly increasing function on $[0,\infty)$. As a result, it implies from (\ref{t4.7}) that
\begin{align}
\int_{Q_R} |u(t,x)|\big(\phi^*_R(t,x)\big)^{\frac{\nu}{\nu+2}}\,d(x,t) &=\int_{Q^*_R} |u(t,x)|\big(\phi^*_R(t,x)\big)^{\frac{\nu}{\nu+2}}\,d(x,t) \nonumber \\
&\le R^{n+2}\,\Phi_q^{-1}\left(\frac{\displaystyle\int_{Q_R} \Phi_q\big( |u(t,x)|\big)\, \big(\phi^*_R(t,x)\big)^{\frac{\nu q_c}{\nu+2}}\,d(x,t)}{R^{n+2}}\right). \label{t4.8}
\end{align}
Collecting the estimates (\ref{t4.1}), (\ref{t4.3}) and (\ref{t4.8}) we obtain
$$ I_R+ \int_{\R^n} u_1(x)\phi_R(0,x)\,dx \lesssim R^n\,\Phi_q^{-1}\left(\frac{\displaystyle\int_{Q_R} \Phi_q\big( |u(t,x)|\big)\, \big(\phi^*_R(t,x)\big)^{\frac{\nu q_c}{\nu+2}}\,d(x,t)}{R^{n+2}}\right). $$
In the same manner one also gets
$$ J_R+ \int_{\R^n} v_1(x)\phi_R(0,x)\,dx \lesssim R^n\,\Phi_p^{-1}\left(\frac{\displaystyle\int_{Q_R} \Phi_p\big( |v(t,x)|\big)\, \big(\phi^*_R(t,x)\big)^{\frac{\nu p_c}{\nu+2}}\,d(x,t)}{R^{n+2}}\right). $$
Thanks to the assumption (\ref{optimaldata}), we have
$$ \int_{\R^n} u_1(x)\phi_R(0,x)\,dx >0 \quad \text{ and }\quad \int_{\R^n} v_1(x)\phi_R(0,x)\,dx >0 $$
for all $R > R_0$, where $R_0$ is a sufficiently large, positive constant. Thus, it holds
\begin{align}
I_R &\lesssim R^n\,\Phi_q^{-1}\left(\frac{\displaystyle\int_{Q_R} \Phi_q\big( |u(t,x)|\big)\, \big(\phi^*_R(t,x)\big)^{\frac{\nu q_c}{\nu+2}}\,d(x,t)}{R^{n+2}}\right), \label{t4.9} \\ 
J_R &\lesssim R^n\,\Phi_p^{-1}\left(\frac{\displaystyle\int_{Q_R} \Phi_p\big( |v(t,x)|\big)\, \big(\phi^*_R(t,x)\big)^{\frac{\nu p_c}{\nu+2}}\,d(x,t)}{R^{n+2}}\right), \label{t4.10}
\end{align}
for all $R > R_0$. Next, for $\lambda>0$ we define the following auxiliary functions:
\begin{align*}
g_p(\lambda) &= \int_{Q_R} \Phi_p\big( |v(t,x)|\big)\, \big(\phi^*_\lambda(t,x)\big)^{\frac{\nu p_c}{\nu+2}}\,d(x,t) \quad \text{ and }\quad G_p(R)= \int_0^R g_p(\lambda)\lambda^{-1}\,d\lambda, \\ 
g_q(\lambda) &= \int_{Q_R} \Phi_q\big( |u(t,x)|\big)\, \big(\phi^*_\lambda(t,x)\big)^{\frac{\nu q_c}{\nu+2}}\,d(x,t) \quad \text{ and }\quad G_q(R)= \int_0^R g_q(\lambda)\lambda^{-1}\,d\lambda.
\end{align*}
Therefore, we can express
\begin{align*}
G_q(R) &= \int_0^R \left(\int_{Q_R} \Phi_q\big( |u(t,x)|\big)\, \big(\phi^*_\lambda(t,x)\big)^{\frac{\nu q_c}{\nu+2}}\,d(x,t)\right) \lambda^{-1}\,d\lambda \\
&= \int_{Q_R} \Phi_q\big( |u(t,x)|\big)\left(\int_0^R \Big(\varphi^*\Big(\frac{t^2+|x|^4}{\lambda^4}\Big)\Big)^{\nu q_c} \lambda^{-1}\,d\lambda\right)\,d(x,t).
\end{align*}
By performing change of variables $\bar{\lambda}= \frac{t^2+|x|^4}{\lambda^4}$, one has
\begin{align*}
G_q(R) &= C \int_{Q_R} \Phi_q\big( |u(t,x)|\big)\left(\int_{\frac{t^2+|x|^4}{R^4}}^\ity \big(\varphi^*(\bar{\lambda})\big)^{\nu q_c}\, \bar{\lambda}^{-1}\,d\bar{\lambda}\right)\,d(x,t) \\
&\le C \int_{Q_R} \Phi_q\big( |u(t,x)|\big)\left(\int_{\frac{1}{2}}^1 \big(\varphi^*(\bar{\lambda})\big)^{\nu q_c}\, \bar{\lambda}^{-1}\,d\bar{\lambda}\right)\,d(x,t) \\
&\hspace{5cm} \big(\text{since}\,\, \text{supp}\,\varphi^* \subset [1/2,1]\big) \\
&\le C \int_{Q_R} \Phi_q\big( |u(t,x)|\big)\max_{\lambda \in (0,R)}\Big(\varphi\Big(\frac{t^2+|x|^4}{\lambda^4}\Big)\Big)^{\nu q_c}\, \left(\int_{\frac{1}{2}}^1 \bar{\lambda}^{-1}\,d\bar{\lambda}\right)\,d(x,t) \\
&\hspace{3cm} \big(\text{since}\,\, \varphi^* \equiv \varphi \,\,\text{in } [1/2,1]\big) \\
&\le C \int_{Q_R} \Phi_q\big( |u(t,x)|\big) \Big(\varphi\Big(\frac{t^2+|x|^4}{R^4}\Big)\Big)^{\nu q_c}\,d(x,t),
\end{align*}
since $\varphi$ is decreasing. An analogous argument results
$$ G_p(R)\le C \int_{Q_R} \Phi_p\big( |u(t,x)|\big) \Big(\varphi\Big(\frac{t^2+|x|^4}{R^4}\Big)\Big)^{\nu p_c}\,d(x,t). $$
Let us now choose $\nu\ge \max\big\{\frac{2}{p_c-1},\frac{2}{q_c-1}\big\}= \frac{2}{p_c-1}$. Thus, it follows that
\begin{align}
G_p(R) &\le C \int_{Q_R} |v(t,x)|^{p_c}\mu_1\big(|u(t,x)|\big) \Big(\varphi\Big(\frac{t^2+ |x|^4}{R^4}\Big)\Big)^{\nu+2}\,d(x,t)= C I_R, \label{t4.11} \\
G_q(R) &\le C \int_{Q_R} |u(t,x)|^{q_c}\mu_2\big(|u(t,x)|\big) \Big(\varphi\Big(\frac{t^2+ |x|^4}{R^4}\Big)\Big)^{\nu+2}\,d(x,t)= C J_R. \label{t4.12}
\end{align}
Furthermore, the following relations hold:
\begin{align}
G_p'(R)R= g_p(R)= \int_{Q_R} \Phi_p\big( |v(t,x)|\big)\, \big(\phi^*_R(t,x)\big)^{\frac{\nu p_c}{\nu+2}}\,d(x,t), \label{t4.13} \\ 
G_q'(R)R= g_q(R)= \int_{Q_R} \Phi_q\big( |u(t,x)|\big)\, \big(\phi^*_R(t,x)\big)^{\frac{\nu q_c}{\nu+2}}\,d(x,t). \label{t4.14}
\end{align}
Combining the estimates from (\ref{t4.9}) to (\ref{t4.14}) gives
\begin{align*}
\frac{G_p(R)}{C} &\le I_R \le C R^n\,\Phi_q^{-1}\Big(\frac{G_q'(R)}{R^{n+1}}\Big), \\
\frac{G_q(R)}{C} &\le J_R \le C R^n\,\Phi_p^{-1}\Big(\frac{G_p'(R)}{R^{n+1}}\Big),
\end{align*}
which are equivalent to
\begin{align*}
\Phi_p\left(\frac{G_q(R)}{C R^n}\right) &\le \frac{G_p'(R)}{R^{n+1}}, \\
\Phi_q\left(\frac{G_p(R)}{C R^n}\right) &\le \frac{G_q'(R)}{R^{n+1}}, \\
\end{align*}
for all $R> R_0$. Then, recalling the definition of the functions $\Phi_p$ and $\Phi_q$ we derive
\begin{align*}
\left(\frac{G_q(R)}{C R^n}\right)^{p_c}\mu_1\left(\frac{G_q(R)}{C R^n}\right) &\le \frac{G_p'(R)}{R^{n+1}}, \\ 
\left(\frac{G_p(R)}{C R^n}\right)^{q_c}\mu_2\left(\frac{G_p(R)}{C R^n}\right) &\le \frac{G_q'(R)}{R^{n+1}},
\end{align*}
for all $R> R_0$. Consequently, it implies
\begin{align*}
\frac{1}{C R^{n(p_c-1)-1}}\mu_1\left(\frac{G_q(R)}{C R^n}\right)\big(G_q(R)\big)^{p_c} &\le G_p'(R), \\ 
\frac{1}{C R^{n(q_c-1)-1}}\mu_2\left(\frac{G_p(R)}{C R^n}\right)\big(G_p(R)\big)^{q_c} &\le G_q'(R),
\end{align*}
for all $R> R_0$. Thank to the increasing property of the functions $\mu_1,\mu_2$ and $G_p=G_p(R)$, $G_q=G_q(R)$, the following inequalities hold:
\begin{align*}
\frac{1}{C R^{n(p_c-1)-1}}\mu_1\left(\frac{G_q(R_0)}{C R^n}\right)\big(G_q(R)\big)^{p_c} &\le G_p'(R), \\
\frac{1}{C R^{n(q_c-1)-1}}\mu_2\left(\frac{G_p(R_0)}{C R^n}\right)\big(G_p(R)\big)^{q_c} &\le G_q'(R),
\end{align*}
hence,
\begin{align*}
\frac{C}{R^{n(p_c-1)-1}}\mu_1\big(C_0 R^{-n}\big)\big(G_q(R)\big)^{p_c} &\le G_p'(R), \\
\frac{C}{R^{n(q_c-1)-1}}\mu_2\big(C_0 R^{-n}\big)\big(G_p(R)\big)^{q_c} &\le G_q'(R),
\end{align*}
for all $R> R_0$, where $C_0:= \frac{1}{C}\min\big\{G_p(R_0),G_q(R_0)\big\}$. For the simplicity, putting $r:= R$ and denoting
$$ \theta_1(r):= \frac{1}{r^{n(p_c-1)-1}}\mu_1\big(C_0 r^{-n}\big), \quad \theta_2(r):= \frac{1}{r^{n(q_c-1)-1}}\mu_2\big(C_0 r^{-n}\big), $$
we obtain the following system of ordinary differential inequalities for $r>R_0$:
\begin{align}
G_p'(r) &\ge C\,\theta_1(r) \big(G_q(r)\big)^{p_c}, \label{trial.1} \\ 
G_q'(r) &\ge C\,\theta_2(r) \big(G_p(r)\big)^{q_c}. \label{trial.2}
\end{align}
Multiplying (\ref{trial.1}) by $G_q'(r)$ and integrating the resulting inequality by parts over $[R_0,r]$ we may arrive at
\begin{align*}
&G_p(r)G_q'(r)- G_p(R_0)G_q'(R_0)- \int_{R_0}^r G_p(\tau)G_q''(\tau)d\tau \\
&\qquad \ge \frac{C}{p_c+1} \theta_1(r) \big(G_q(r)\big)^{p_c+1}- \frac{C}{p_c+1} \theta_1(R_0) \big(G_q(R_0)\big)^{p_c+1} - \frac{C}{p_c+1} \int_{R_0}^r \theta'_1(\tau) \big(G_q(\tau)\big)^{p_c+1}d\tau.
\end{align*}
As a consequence, it implies that
$$ G_p(r)G_q'(r) \ge C \theta_1(r) \big(G_q(r)\big)^{p_c+1}, $$
that is,
\begin{equation}
G_p(r) \ge \frac{C \theta_1(r) \big(G_q(r)\big)^{p_c+1}}{G_q'(r)} \label{trial.3}
\end{equation}
$$  $$
for $r>R_0$. By plugging (\ref{trial.3}) into (\ref{trial.2}), one gets
$$ G_q'(r) \ge \frac{C\theta_2(r)\, \big(\theta_1(r)\big)^{q_c} \big(G_q(r)\big)^{q_c(p_c+1)}}{\big(G_q'(r)\big)^{q_c}}, $$
which is equivalent to
\begin{align*}
G_q'(r) &\ge C \big(\theta_2(r)\big)^{\frac{1}{q_c+1}}\, \big(\theta_1(r)\big)^{\frac{q_c}{q_c+1}} \big(G_q(r)\big)^{\frac{q_c(p_c+1)}{q_c+1}} \\ 
&\quad = \frac{C}{r}\big(\mu_1(C_0 r^{-n})\big)^{\frac{q_c}{q_c+1}}\big(\mu_2(C_0 r^{-n})\big)^{\frac{1}{q_c+1}} \big(G_q(r)\big)^{\frac{q_c(p_c+1)}{q_c+1}}.
\end{align*}
Thus, it follows that
$$ \frac{C}{r}\big(\mu_1(C_0 r^{-n})\big)^{\frac{q_c}{q_c+1}}\big(\mu_2(C_0 r^{-n})\big)^{\frac{1}{q_c+1}} \le \frac{G_q'(r)}{\big(G_q(r)\big)^{\frac{q_c(p_c+1)}{q_c+1}}}. $$
Integrating two sides of the last estimate over $[R_0,R]$ leads to
\begin{align*}
C \int_{R_0}^R \frac{1}{r} \big(\mu_1(C_0 r^{-n})\big)^{\frac{q_c}{q_c+1}}\big(\mu_2(C_0 r^{-n})\big)^{\frac{1}{q_c+1}}\,dr &\le \int_{R_0}^R \frac{G_q'(r)}{\big(G_q(r)\big)^{\frac{q_c(p_c+1)}{q_c+1}}}\,dr \\
&\quad= -\frac{p_cq_c-1}{q_c+1} \big(G_q(r)\big)^{-\frac{p_cq_c-1}{q_c+1}}\Big|_{r=R_0}^{r=R} \le \frac{2}{n}\big(G_q(R_0)\big)^{-\frac{2}{n}},
\end{align*}
where we notice that $\frac{q_c+1}{p_cq_c-1}= \frac{n}{2}$. For this reason, we pass $R \to \ity$ to derive
$$ C \int_{R_0}^\ity \frac{1}{r} \big(\mu_1(C_0 r^{-n})\big)^{\frac{q_c}{q_c+1}}\big(\mu_2(C_0 r^{-n})\big)^{\frac{1}{q_c+1}}\,dr \le \frac{2}{n}\big(G_q(R_0)\big)^{-\frac{2}{n}}. $$
Finally, carrying out change of variables $s= C_0 r^{-n}$ gives
$$ C\int_0^{C_0 R_0^{-n}} \frac{1}{s} \big(\mu_1(s)\big)^{\frac{q_c}{q_c+1}}\big(\mu_2(s)\big)^{\frac{1}{q_c+1}}\,ds \le \frac{2}{n}\big(G_q(R_0)\big)^{-\frac{2}{n}}. $$
This contradicts to the assumption (\ref{blowupCon}). Summarizing, the proof of Theorem \ref{dloptimal1} is completed.
\end{proof}

%=================================================================={References}


\begin{thebibliography}{00}

\bibitem{EbertGirardiReissig} M.R. Ebert, G. Girardi, M. Reissig, \textit{Critical regularity of nonlinearities in semilinear classical damped wave equations}, Math. Ann., (2019), https://doi.org/10.1007/s00208-019-01921-5.
\bibitem{ReissigEbert} M.R. Ebert, M. Reissig, ``Methods for partial differential equations, qualitative properties of solutions, phase space analysis, semilinear models'', Birkh\"auser, 2018.
\bibitem{Friedman} A. Friedman, ``Partial differential equations'', Krieger, New York, 1976.
\bibitem{Matsumura} A. Matsumura, \textit{On the asymptotic behavior of solutions of semi-linear wave equations}, Publ. Res. Inst. Math. Sci., 12 (1976), 169-189.
\bibitem{Narazaki} T. Narazaki, \textit{Global solutions to the Cauchy problem for the weakly coupled system of damped wave equations}, Discrete Contin. Dyn. Syst. Suppl., (2009), 592-601.
\bibitem{NishiharaWakasugi} K. Nishihara, Y. Wakasugi, \textit{Critical exponent for the Cauchy problem to the weakly coupled damped wave systems}, Nonlinear Anal., 108 (2014), 249-259.
\bibitem{SunWang} F. Sun, M. Wang, \textit{Existence and nonexistence of global solutions for a non-linear hyperbolic system with damping}, Nonlinear Anal., 66 (2007), 2889-2910.
\bibitem{Takeda} H. Takeda, \textit{Global existence and nonexistence of solutions for a system of nonlinear damped wave equations}, J. Math. Anal. Appl., 360 (2009), 631-650.
\bibitem{TodorovaYordanov} G. Todorova, B. Yordanov, \textit{Critical exponent for a nonlinear wave equation with damping}, J. Differential Equations, 174 (2001), 464-489.
\bibitem{Zhang} Q.S. Zhang, \textit{A blow-up result for a nonlinear wave equation with damping: the critical case}, C. R. Acad. Sci. Paris S\'{e}r. I Math., 333 (2001), 109-114.

\end{thebibliography}
\end{document}